\RequirePackage{ifpdf}
\ifpdf 
\documentclass[pdftex]{sigma}
\else
\documentclass{sigma}
\fi

\newtheorem{prop}[theorem]{Proposition}

\newtheorem{lem}[theorem]{Lemma}
\newtheorem{def-lem}[theorem]{Definition-Lemma}
\newtheorem{def-prop}[theorem]{Definition-Proposition}

\newtheorem{cor}[theorem]{Corollary}

\def\a{\alpha}\def\Ar{{\rm A}}\def\Ab{\bar{\rm A}}\def\b{\beta}\def\bb{\mathfrak{b}}
\def\Cprime{C'}\def\ds{\displaystyle}\newcommand{\fracds}[2]{{\ds \frac{#1}{#2}}}
\def\f{\mathfrak{k}}\def\h{\mathfrak{h}}\def\hs{\mathring{z}}\def\g{\mathfrak{g}}\def\sl{\mathfrak{sl}}\def\gl{\mathfrak{gl}}
\def\Ib{\tts\overline{\nns\rm I}_+\tts}\def\J{\mathrm{J}}\def\Jb{\,\overline{\!\rm I}_-\ts}\def\k{\mathfrak{k}}
\def\lcd{\ts,\ldots,}\def\mult{\diamond}\def\n{\mathfrak{n}}\def\nminus{\mathfrak{n}_-}
\def\nplus{\mathfrak{n}}
\def\nns{\hskip-.5pt}\def\p{\mathfrak{p}}
\def\pp{\widetilde{\phantom{\,}p\phantom{\,}\phantom{\ }}\phantom{\!}\phantom{\!}\phantom{\!}}
\def\q{\check{\mathrm q}}\def\Q{\mathrm{Q}}\def\r#1{(\ref{#1})}\let\rf=\r
\def\s{z}\def\S{\mathrm {S}}\def\si{\sigma}\def\Sp{\S^{\nplus}}\def\sy{\acute{\sigma}}\def\th{\mathring{h}}
\def\tphi{A}\def\tpsi{A'}\def\ttphi{B}\def\ttpsi{B'}\def\ts{\hskip1pt}\def\tts{\hskip.5pt}
\def\tt{\mathring{t}}\def\U{\mathrm{U}}\def\Uh{\overline{\U}(\h)}\def\UU{\overline{\U}}\def\UUh{{\cal{D}}(\h)}
\def\UUb{\UU^{12}(\bb)}\def\V{\mathrm {V}}\def\ve{\varepsilon}\def\Z{\mathrm{Z}}
\def\Zp{\Z^{\nplus}}\def\z{z}\let\mmult=\mult
\newcommand{\ZZ}{{\mathbb Z}}\def\Hom{\operatorname{Hom}}\def\I{\mathrm{I}}

\numberwithin{equation}{section}

\begin{document}

\allowdisplaybreaks

\renewcommand{\PaperNumber}{064}

\FirstPageHeading

\ShortArticleName{Structure Constants of Diagonal Reduction Algebras  of $\gl$ Type}

\ArticleName{Structure Constants of Diagonal Reduction\\ Algebras  of $\gl$ Type}

\Author{Sergei KHOROSHKIN~$^{a,b}$ and  Oleg OGIEVETSKY~$^{c,d,e}$}

\AuthorNameForHeading{S.~Khoroshkin and  O.~Ogievetsky}

\Address{$^{a)}$~Institute of Theoretical and Experimental Physics, 117218 Moscow, Russia}
\EmailDD{\href{mailto:khor@itep.ru}{khor@itep.ru}}

\Address{$^{b)}$~Higher School of Economics,  20 Myasnitskaya Str., 101000 Moscow, Russia}

\Address{$^{c)}$~J.-V. Poncelet French-Russian Laboratory, UMI 2615 du CNRS,\\
\hphantom{$^{c)}$}~Independent University of Moscow, 11 B.~Vlasievski per.,
119002 Moscow, Russia}
\EmailDD{\href{mailto:oleg.ogievetsky@gmail.com}{oleg.ogievetsky@gmail.com}}

\Address{$^{d)}$~Centre de Physique Th\'eorique\footnote{Unit\'e Mixte de Recherche
(UMR 6207) du CNRS et des Universit\'es Aix--Marseille I,
Aix--Marseille II et du Sud Toulon -- Var; Laboratoire Affili\'e \`a la FRUMAM (FR 2291)},
Luminy, 13288 Marseille, France}

\Address{$^{e)}$~On leave of absence
from P.N.~Lebedev Physical Institute, Theoretical Department,\\
\hphantom{$^{e)}$}~53 Leninsky Prospekt, 119991 Moscow, Russia}

\ArticleDates{Received January 14, 2011, in f\/inal form June 27, 2011;  Published online July 09, 2011}

\Abstract{We describe, in terms of generators and relations, the reduction algebra, related to the diagonal embedding of the Lie algebra
$\gl_n$ into $\gl_n\oplus\gl_n$. Its representation theory is related to the theory of decompositions of tensor products of $\gl_n$-modules.}


\Classification{16S30; 17B35}

{\footnotesize {\tableofcontents}}

\section{Introduction}

This paper completes the work \cite{KO3}: it contains a derivation of basic relations for the diagonal reduction algebras of $\gl$ type, their
low dimensional examples and properties.

Let $\g$ be a Lie algebra, $\f\subset \g$ its reductive Lie subalgebra and $V$ an irreducible f\/inite-dimensional $\g$-module,
which decomposes, as an $\f$-module, into a direct sum of irreducible $\k$-modules $V_i$ with certain multi\-plicities $m_i$,
\begin{gather}
 \label{intro0}V\approx \sum_iV_i\otimes W_i .
\end{gather}
Here  $W_i=\Hom_\f(V_i,V)$ are the spaces of multiplicities, $m_i=\dim W_i$. While the multiplici\-ties~$m_i$ present certain combinatorial data, the spaces
$W_i$ of multiplicities itself may exhibit a~`hidden structure'  of modules over certain special algebras~\cite{C}. The well-known example is the Olshanski
\textit{centralizer construction}~\cite{Ol}, where  $\g=\gl_{n+m}$, $\f=\gl_m$ and the spaces $W_i$ carry the structure of irreducible Yangian $Y(\gl_n)$-modules.

In general, the multiplicity spaces $W_i$ are irreducible modules over the centralizer $\U(\g)^\f$ of $\f$ in the universal enveloping algebra
$\U(\g)$~\cite{LM}. However, these centralizers have a rather complicated algebraic structure and are hardly convenient for applications. Besides, under
the above assumptions, the direct sum $W=\oplus_i W_i$ becomes  a module over the {\it reduction} (or Mickelsson) algebra. The reduction algebra is def\/ined
as follows. Suppose $\k$ is given with  a~triangular decomposition
\begin{gather}
\label{intro1} \k=\nminus+\h+\nplus  .
\end{gather}
Denote by $\I_+$ the left ideal of $\Ar :=\U(\g)$, generated by elements of $\nplus$, $\I_+ :=\Ar\nplus$ . Then the reduction algebra $\Sp(\Ar)$, related to the pair
$(\g,\k)$, is def\/ined as the quotient $\mathrm{Norm}(\I_+)/\I_+$ of the normalizer of the ideal $\I_+$ over $\I_+$. It is equipped with a~natural structure of the
associative algebra.  By def\/inition, for any $\g$-module $V$ the space $V^\nplus$ of vectors, annihilated by $\nplus$, is a~module over~$\Sp(\Ar)$. Since $V$ is
f\/inite-dimensional, $V^\nplus$ is isomorphic to $\oplus_iW_i$, so the latter space can be viewed as an~$\Sp(\Ar)$-module as well; the zero-weight component of~$\Sp(\Ar)$, which contains a quotient of the centralizer~$\U(\g)^\f$, preserves each multiplicity space~$W_i$. The representation theory of the reduction
algebra~$\Sp(\Ar)$ is closely related to the theory of branching rules $\g\downarrow \k$ for the restrictions of representations of~$\g$ to~$\k$.

The reduction algebra simplif\/ies  after the localization over the multiplicative set generated by elements $h_\gamma+k$, where $\gamma$
ranges through the set of roots of $\k$, $k\in\ZZ$, and $h_\gamma$ is the coroot corresponding to $\gamma$. Let $\Uh$ be the localization of the universal
enveloping algebra $\U(\h)$ of the Cartan subalgebra $\h$ of $\k$ over the above multiplicative set. The localized reduction algebra~$\Zp(\Ar)$ is
an algebra over the commutative ring~$\Uh$; the principal part of the def\/ining relations is quadratic but the relations may contain linear or degree 0 terms, see~\cite{Zh,KO}.

Besides, the reduction algebra admits another description as a (localized) double coset space $\nminus\Ar\backslash \Ar/\Ar \nplus$ endowed with the
multiplication map def\/ined by means of the insertion of the extremal projector~\cite{KO} of Asherova--Smirnov--Tolstoy~\cite{AST}. The centralizer $\Ar^\f$ is a subalgebra of~$\Zp(\Ar)$.

It was shown in \cite{KO3} that the general reduction algebra $\Zp(\Ar)$ admits a presentation over $\Uh$ such that the def\/ining relations are
ordering relations for the generators, in an arbitrary order, compatible with the natural partial order on~$\h^*$.
The set of ordering relations for the general reduction algebra $\Zp(\Ar)$ was shown in \cite{KO3} to be ``algorithmically ef\/f\/icient'' in the sense that any
expression in the algebra can be ordered with the help of this set.

The structure constants of the reduction algebra are in principle determined with the help of the extremal projector $P$ or the tensor $\J$ studied by
Arnaudon,  Buf\/fenoir, Ragoucy and Roche~\cite{ABRR}. However the explicit description of the algebra is hardly achievable directly.

In the present paper, we are interested in the special restriction problem, when $\g$ is the direct sum of two copies of a
reductive Lie algebra ${\mathfrak a}$ and $\f$ is the diagonally embedded ${\mathfrak a}$. The resulting reduction algebra for the symmetric pair
$({\mathfrak a}\oplus {\mathfrak a}, {\mathfrak a})$ we call {\em diagonal reduction} algebra~${\mathrm{DR}}({\mathfrak a})$ of~${\mathfrak a}$.  The theory of branching rules for ${\mathfrak a}\oplus {\mathfrak a}\downarrow {\mathfrak a}$ is the theory of
decompositions of the tensor products of ${\mathfrak a}$-modules into a direct sum of irreducible ${\mathfrak a}$-modules.

We restrict ourselves here to the Lie algebra ${\mathfrak a}=\gl_n$ of the general linear group. In this situation f\/inite-dimensional irreducible modules over $\g$
are tensor products of two irreducible $\gl_n$-modules, the decomposition \rf{intro0} is the decomposition of the tensor product into the direct sum of irreducible
$\gl_n$-modules, and the multiplicities $m_i$ are the Littlewood--Richardson coef\/f\/icients.

The reduction algebra ${\mathrm{DR}}(\gl_n)$ for brevity will be denoted further by $\Z_n$.

In \cite{KO3} we suggested a set $\mathfrak{R}$ of relations for the algebra $\Z_n$.  We demonstrated that the set $\mathfrak{R}$ is equivalent, over $\Uh$,
to the set of the def\/ining ordering relations provided that all relations from the set $\mathfrak{R}$ are valid.

The main goal of the present paper is the verif\/ication of all relations from the system $\mathfrak{R}$. There are two
principal tools in our derivation. First, we use the braid group action by the  Zhelobenko automorphisms of reduction algebras \cite{Zh,KO}. Second, we
employ the stabilization phenomenon, discovered in \cite{KO3},
for the multiplication rule and for the def\/ining relations with respect to the standard embeddings
$\gl_n\hookrightarrow \gl_{n+1}$;  stabilization provides a natural way of extending relations for $\Z_{n}$ to relations for $\Z_{n+1}$ ($\Z_{n}$ is not a
subalgebra of $\Z_{n+1}$). We perform necessary calculations for low $n$ (at most $n=4$); the braid group action and the stabilization law allow
to extend the results for general $n$.

As an illustration, we write down the complete lists of def\/ining relations in the form of ordering relations for the reduction algebras ${\mathrm{DR}}(\sl_3)$
and ${\mathrm{DR}}(\sl_2)$. Although for a f\/inite $n$ the task of deriving the set of def\/ining (ordering) relations for ${\mathrm{DR}}(\sl_n)$ is achievable in a f\/inite time, it is useful to have the list of relations for small $n$ in front of the eyes.

We return to the stabilization and cut phenomena and make more precise statements concerning now the embedding of the Lie algebra $\gl_n\oplus\gl_1$
into the Lie algebra $\gl_{n+1}$ (more generally, of $\gl_n\oplus\gl_m$ into $\gl_{n+m}$). As a consequence we f\/ind that cutting preserves the centrality:
the cut of a central element of the algebra $\Z_{n+m}$ is central in the algebra $\Z_n\otimes\Z_m$.
We also show that, similarly to the Harish-Chandra map, the restriction of the cutting to the center is a homomorphism. As an example, we
derive the Casimir operators for the algebra ${\mathrm{DR}}(\sl_2)$ by cutting the Casimir operators for the algebra ${\mathrm{DR}}(\sl_3)$.

The relations in the diagonal reduction algebra have a quadratic and a degree zero part. The algebra, def\/ined by the homogeneous quadratic part
of the relations, tends, in a quite simple regime, to a commutative algebra (the homogeneous algebra can be thus considered as a ``dynamical''
deformation of a commutative algebra; ``dynamical'' here means that the left and right multiplications by elements of the ring $\U (\h )$ dif\/fer). This observation
about the limit is used in the proof in~\cite{KO3} of the completeness of the set of derived relations over the f\/ield of
fractions of $\U(\h)$. We prove the completeness by establishing the equivalence between the set of derived relations and the set of ordering relations.

The stabilization law enables one to give a def\/inition of the reduction ``algebra'' $\Z_\infty$ related to the diagonal embedding of the inductive
limit $\gl_\infty$ of $\gl_n$ into $\gl_\infty\oplus \gl_\infty$ (strictly speaking, $\Z_\infty$ is not an algebra, some relations have an inf\/inite number of terms).

We also discuss the diagonal reduction algebra for the special linear Lie algebra $\sl_n$; it is a~direct tensor factor in~$\Z_{n}$.

Such a precise description, as the one we give for~$\Z_{n}$, is known for a few examples of the reduction algebras: the most known is related to the
embedding of $\gl_n$ to $\gl_{n+1}$ \cite{Zh}. Its representation theory was used for the derivation of precise formulas for the action of the generators
of $\gl_n$ on the Gelfand--Zetlin basic vectors~\cite{AST2}. The reduction algebra for the pair $(\gl_n,\gl_{n+1})$ is based on the root embedding
$\gl_n\subset\gl_{n+1}$ of Lie algebras. In contrast to this example, the  diagonal reduction algebra ${\mathrm{DR}}({\mathfrak a})$ is
based on the diagonal embedding of ${\mathfrak a}$ into
${\mathfrak a}\oplus {\mathfrak a}$, which is not a root embedding of reductive Lie algebras.

\section{Notation}\label{section-notation}

Let ${\mathcal{E}}_{ij}$, $i,j=1\lcd n$, be the standard generators of the Lie algebra $\gl_n$, with the commutation relations
\begin{gather*}
[{\mathcal{E}}_{ij},{\mathcal{E}}_{kl}]=\delta_{jk}{\mathcal{E}}_{il}-\delta_{il}{\mathcal{E}}_{kj} ,
\end{gather*}
where $\delta_{jk}$ is the Kronecker symbol. We shall also use the root notation ${\mathcal{H}}_\alpha$, ${\mathcal{E}}_{\alpha}$, ${\mathcal{E}}_{-\alpha}$, \dots\ for elements of~$\gl_n$.

Let ${\mathcal{E}}_{ij}^{(1)}$ and ${\mathcal{E}}_{ij}^{(2)}$, $i,j=1\lcd n$, be the standard generators of the two copies of the Lie algebra $\gl_n$ in $\g:=\gl_n\oplus\gl_n$,
\begin{gather*}
[{\mathcal{E}}_{ij}^{(a)},{\mathcal{E}}_{kl}^{(b)}]=\delta_{ab}\big(\delta_{jk}{\mathcal{E}}_{il}^{(a)}-\delta_{il}{\mathcal{E}}_{kj}^{(a)}\big) .\end{gather*}
Set
\begin{gather*}
 e_{ij}:={\mathcal{E}}_{ij}^{(1)}+{\mathcal{E}}_{ij}^{(2)} ,\qquad E_{ij}:={\mathcal{E}}_{ij}^{(1)}-{\mathcal{E}}_{ij}^{(2)}
 .
 \end{gather*}
The elements $e_{ij}$ span the diagonally embedded Lie algebra $\k\simeq\gl_n$, while $E_{ij}$ form an adjoint $\k$-module $\mathfrak{p}$. The
Lie algebra $\k$ and the space $\mathfrak{p}$ constitute a symmetric pair, that is, $[\k,\k]\subset \k$, $[\k,\p]\subset \p$, and $[\p,\p]\subset \k$:
\begin{gather*}
 [e_{ij},e_{kl}] = \delta_{jk}e_{il}-\delta_{il}e_{kj}  ,\qquad
  [e_{ij},E_{kl}] = \delta_{jk}E_{il}-\delta_{il}E_{kj}  ,\qquad
   [E_{ij},E_{kl}] = \delta_{jk}e_{il}-\delta_{il}e_{kj} .
\end{gather*}
In the sequel, $h_a$ means the element $e_{aa}$ of the Cartan subalgebra $\h$ of the subalgebra $\k\in\gl_n\oplus\gl_n$ and $h_{ab}$ the
element $e_{aa}-e_{bb}$.

Let $\{\ve_a\}$ be the basis of $\h^*$ dual to the basis $\{ h_a\}$ of $\h$, $\ve_a(h_b)=\delta_{ab}$. We shall use as well the root notation $h_\alpha$,
$e_{\alpha}$, $e_{-\alpha}$ for elements of $\k$, and $H_\alpha$, $E_{\alpha}$, $E_{-\alpha}$ for elements of $\p$.

The Lie subalgebra $\nplus$ in the triangular decomposition \rf{intro1} is spanned by the root vectors~$e_{ij}$ with $i<j$ and the Lie subalgebra $\nminus$ by the
root vectors $e_{ij}$ with $i>j$. Let $\bb_+$ and $\bb_-$
be the corresponding Borel subalgebras, $\bb_+=\h\oplus\nplus$ and $\bb_-=\h\oplus\nminus$. Denote by $\Delta_+$ and $\Delta_-$ the sets of positive and
negative roots in the root system $\Delta=\Delta_+\cup \Delta_-$ of $\k$: $\Delta_+$ consists of roots $\ve_i-\ve_j$ with $i<j$ and $\Delta_-$ consists of roots
$\ve_i-\ve_j$ with $i>j$.
Let $\Q$ be the root lattice, $\Q:=\{\gamma\in\h^*\,|\,\gamma=\sum_{\alpha\in\Delta_+, n_\a\in\ZZ}n_\a \a\}$. It contains the
positive cone $\Q_+$,
\begin{gather*}
\Q_+:=\bigg\{\gamma\in\h^*\,|\,\gamma= \sum_{\alpha\in\Delta_+, n_\a\in\ZZ, \; n_\a\geq 0}n_\a \a\bigg\} .
\end{gather*}
For $\lambda,\mu\in\h^*$, the notation
\begin{gather}
\lambda>\mu\label{paor}
\end{gather}
means that the dif\/ference $\lambda-\mu$ belongs to~$\Q_+$, $\lambda-\mu\in\Q_+$. This is a partial order in~$\h^*$.

We f\/ix the following action of the cover of the symmetric group $\S_n$ (the Weyl group of the
diagonal $\k$) on the Lie algebra $\gl_n\oplus\gl_n$ by automorphisms
\begin{gather*}
\sy_i(x):=\mathrm{Ad}_{\exp(e_{i,i+1})}\mathrm{Ad}_{\exp(-e_{i+1,i})}\mathrm{Ad}_{\exp(e_{i,i+1})}(x) ,
\end{gather*}
so that
\begin{gather*}\sy_i(e_{kl})=(-1)^{\delta_{ik}+\delta_{il}}e_{\sigma_i(k)\sigma_i(l)} ,
\qquad\sy_i(E_{kl})=(-1)^{\delta_{ik}+\delta_{il}}E_{\sigma_i(k)\sigma_i(l)}.
\end{gather*}
Here $\sigma_i=(i,i+1)$ is an elementary transposition in the symmetric group. We extend naturally the above action of the cover of $\S_n$ to the action by
automorphisms on the associative algebra $\Ar\equiv\Ar_n:=\U(\gl_n)\otimes \U(\gl_n)$. The restriction of this action to $\h$ coincides
with the natural action $\si(h_{k})=h_{\si(k)}$, $\si\in\S_n$,  of the Weyl group on the Cartan subalgebra.

Besides, we use the shifted action of $\S_n$ on the polynomial algebra $\U(\h)$ (and its localizations) by automorphisms; the shifted action is def\/ined by
\begin{gather}
\label{not1}
\si\circ h_k:= h_{\si(k)}+k-\si(k) ,\qquad k=1,\dots ,n;\qquad \si\in\S_n .
\end{gather}
It becomes the usual action for the variables
\begin{gather}
 \th_k:=h_k-k ,\qquad \th_{ij}:=\th_i-\th_j; \label{hrond}
\end{gather}
by \rf{not1} for any $\si\in\S_n$ we have
\begin{gather*}
\si\circ\th_k=\th_{\si(k)} ,\qquad \si\circ\th_{ij}=\th_{\si(i)\si(j)} .
\end{gather*}

It will be sometimes convenient to denote the commutator $[a,b]$ of two elements $a$ and $b$ of an associative algebra by
\begin{gather}
 \hat{a}(b):=[a,b] .\label{dehaco}
\end{gather}

\section[Reduction algebra $\Z_n$]{Reduction algebra $\boldsymbol{\Z_n}$}\label{section2}

In this section we recall the def\/inition of the reduction algebras, in particular the diagonal reduction algebras of the $\gl$ type.
We introduce the order for which
the ordering relations for the algebra $\Z_n$ will be discussed.
The formulas for the Zhelobenko automorphisms for the algebra $\Z_n$ are given; some basic facts about the standard involution, anti-involution
and central elements for the algebra $\Z_n$ are presented at the end of the section.

{\bf 1.} Let $\Uh$ and $\Ab$ be the rings of fractions of the algebras $\U(\h)$ and $\Ar$ with respect to the multiplicative set, generated by elements
\begin{gather*}
 h_{ij}+l ,\qquad l\in\ZZ  , \quad 1\leq i<j\leq n  .
 \end{gather*}

Def\/ine $\Z_n$ to be the double coset space of $\Ab$ by its left ideal $\Ib:=\Ab\nplus$, generated by elements of $\nplus$, and the right ideal
$\Jb:=\nminus\Ab$, generated by elements of $\nminus$, $\Z_n:=\Ab/(\Ib+\Jb)$.

The space $\Z_n$ is an associative algebra with respect to the multiplication map
\begin{gather}
\label{not5a}
a\mult b:=aP b  .
\end{gather}
Here $P$ is the extremal projector \cite{AST} for the diagonal $\gl_n$.
It is an element of a certain extension of the algebra $\U(\gl_n)$ satisfying the relations $e_{ij}P=Pe_{ji}=0$ for all $i$ and
$j$ such that $1\leq i<j\leq n$.

The algebra  $\Z_n$ is a particular example of a reduction algebra; in our context, $\Z_n$ is
def\/ined by the coproduct (the diagonal inclusion) $\U(\gl_n)\to\Ar$.

{\bf 2.}  The main structure theorems for the reduction algebras are given in \cite[Section~2]{KO3}.

In the sequel we choose a weight linear basis $\{p_K\}$ of $\p$ ($\p$ is the $\k$-invariant complement to~$\k$ in~$\g$,
$\g =\k +\p$) and equip it with a total order $\prec$. The total order~$\prec$ will be compatible with the partial order $<$ on
$\h^*$, see~\rf{paor}, in the sense that $\mu_K<\mu_L\ \Rightarrow \ p_K\prec p_L$. We shall sometimes write
$I\prec J$ instead of $p_I\prec p_J$. For an arbitrary element $a\in\Ab$ let $\widetilde{a}$ be its image in the reduction
algebra; in particular, $\pp_K$ is the image in the reduction algebra of the basic vector $p_K\in\p$.

{\bf 3.} In our situation we choose the set of vectors $E_{ij}$, $i,j=1,\dots ,n$, as a basis of the space~$\p$. The weight of $E_{ij}$ is $\ve_i-\ve_j$.
The compatibility of a total order $\prec$ with the partial order~$<$ on~$\h^*$ means the condition
\begin{gather*}
E_{ij}\prec E_{kl}\qquad \text{if}\quad i-j>k-l .
\end{gather*}
The order in each subset $\{E_{ij}|i-j=a\}$ with a f\/ixed $a$ can be chosen arbitrarily. For instance, we can set
\begin{gather}\label{not4}
E_{ij}\prec E_{kl}\qquad \text{if}\quad i-j>k-l\quad\text{or}\quad i-j=k-l\quad\text{and}\quad i>k .
\end{gather}

Denote the images of the elements $E_{ij}$ in $\Z_n$ by $\s_{ij}$. We use also the notation $t_i$ for the elements $\s_{ii}$ and $t_{ij}:=t_i-t_j$ for the elements
$\s_{ii}-\s_{jj}$. The order \rf{not4} induces as well the order on the generators $\s_{ij}$ of the algebra $\Z_n$:
\begin{gather*} 
\s_{ij}\prec\s_{kl}\ \Leftrightarrow\ E_{ij}\prec E_{kl}  .
\end{gather*}
The  statement (d) in the paper \cite[Section~2]{KO3} implies an existence of structure constants $\mathrm{B}_{(ab),(cd),(ij),(kl)} \in \Uh$ and $\mathrm{D}_{(ab),(cd)}\in\Uh$
such that for any $a,b,c,d=1,\ldots,n$ we have
\begin{gather}\label{not6}
\s_{ab}\mult\s_{cd}=\sum_{i,j,k,l:\s_{ij}\preceq\s_{kl}}\mathrm{B}_{(ab),(cd),(ij),(kl)}\s_{ij}\mult\s_{kl}+\mathrm{D}_{(ab),(cd)}  .
\end{gather}
In particular, the algebra $\Z_n$ (in general, the reduction algebra related to a symmetric pair $(\k,\p)$, $\g :=\k +\p$) is $\ZZ_2$-graded; the degree of $\s_{ab}$ is~1 and the degree of any element from $\Uh$ is~0.

The relations \rf{not6} together with the weight conditions
\begin{gather*}
 [h,\s_{ab}]=(\ve_a-\ve_b)(h) \s_{ab}
\end{gather*}
are the def\/ining relations for the algebra $\Z_n$.

Note that the denominators of the structure constants $\mathrm{B}_{(ab),(cd),(ij),(kl)}$ and
 $\mathrm{D}_{(ab),(cd)}$ are pro\-ducts
of linear factors of the form $\th_{ij}+\ell$, $i<j$, where $\ell\geq -1$ is an integer, see~\cite{KO3}.

{\bf 4.} The algebra $\Z_n$ can be equipped with the action of Zhelobenko automorphisms \cite{KO}.  Denote by $\q_{i}$  the Zhelobenko
automorphism $\q_i:\Z_n\to\Z_n$  corresponding to the transposition $\si_i\in\S_n$. It is def\/ined as follows~\cite{KO}. First we def\/ine a map
$\q_i:\Ar\to \Ab/\Ib$ by
\begin{gather} \label{not7}
\q_i(x):=\sum_{k\geq 0}\frac{(-1)^k}{k!} \hat{e}_{i,i+1}^k(\sy_i(x)) e_{i+1,i}^k   \prod_{a=1}^k(h_{i,i+1}-a+1)^{-1}
\quad \mod\Ib  .
\end{gather}
Here $\hat{e}_{i,i+1}$ stands for the adjoint action of the element $e_{i,i+1}$, see~\rf{dehaco}.
The operator $\q_i$ has the property
\begin{gather}
\label{not2}
\q_i(hx)=(\si_i\circ h)\q_i(x)
\end{gather}
for any $x\in\Ar$ and  $h\in\h$; $\si\circ h$ is def\/ined in~\rf{not1}. With the help of~(\ref{not2}), the map $\q_i$ can be extended
to the map (denoted by the same symbol) $\q_i:\Ab\to \Ab/\Jb$ by
setting $\q_i(a(h)x)=(\si_i\circ a(h))\q_i(x)$ for any $x\in\Ar$ and
$a(h)\in\Uh$. One can further prove that $\q_i(\Ib)=0$ and $\q_i(\Jb)\subset (\Jb+\Ib)/\Ib$,
so that $\q_i$ can be viewed as a linear operator $\q_i:\Z_n\to\Z_n$. Due to~\cite{KO},
this is an algebra automorphism, satisfying~\rf{not2}.

The operators $\q_i$ satisfy the braid group relations~\cite{Zh}:
\begin{gather*}
\q_i\q_{i+1}\q_i=\q_{i+1}\q_{i}\q_{i+1}  ,\\
\q_i\q_j=\q_j\q_i  ,\quad |i-j|>1  ,
\end{gather*}
and the inversion relation \cite{KO}:
\begin{gather}
\label{invr}
\q_i^2(x)=\frac{1}{h_{i,i+1}+1}  \sy_i^2(x)  (h_{i,i+1}+1)  ,\qquad x\in\Z_n  .
\end{gather}
In particular, $\q_i^2(x)=x$ if $x$ is of zero weight.

{\bf 5.} The Chevalley anti-involution $\epsilon$ in $\U(\gl_n\oplus\gl_n)$,
$\epsilon(e_{ij}):=e_{ji}$, $\epsilon(E_{ij}):=E_{ji}$, induces the anti-involution $\epsilon$ in
the algebra $\Z_n$:
\begin{gather}
\epsilon(\s_{ij})=\s_{ji}  ,\qquad \epsilon(h_k)=h_k  .\label{anep}
\end{gather}
Besides, the outer automorphism of the Dynkin diagram of $\gl_n$ induces the involutive automorphism
$\omega$ of $\Z_n$,
\begin{gather}
\label{not2a}
\omega(\s_{ij}) =(-1)^{i+j+1}\s_{j'i'}  ,\qquad \omega (h_k)=-h_{k'}  ,
\end{gather}
where $i'=n+1-i$. The operations $\epsilon$ and $\omega$ commute, $\epsilon\omega =\omega\epsilon$.

Central elements of the subalgebra $\U(\gl_n)\otimes 1\subset \Ar$, generated by $n$ Casimir
operators of degrees $1\lcd n$,  as well as central elements of the subalgebra
$1\otimes \U(\gl_n)\subset \Ar$ project to central elements of the algebra $\Z_n$. In particular,
central elements of degree $1$ project to central elements
\begin{gather}
I^{(n,h)}:=h_1+\dots +h_n\label{clih}
\end{gather}
and
\begin{gather}
\label{clit} I^{(n,t)}:=t_1+\dots +t_n
\end{gather}
of the algebra $\Z_n$. The dif\/ference of central elements of degree two projects to the central element
\begin{gather}
\label{drclit}
\sum_{i=1}^n(h_i-2i)t_i
\end{gather}
of the algebra $\Z_n$. The images of other Casimir operators are more complicated.

\section{Main results}

This section contains the principal results of the paper. We f\/irst give
preliminary information on the new basis in which the def\/ining relations for the
algebra $\Z_n$ can be written down in an economical fashion. The braid group action on the new generators is then
explicitly given in Subsection~\ref{brgac}. The complete set of the def\/ining relations for the algebra $\Z_n$ is written down in Subsection~\ref{section3.3}. The regime for which both
the set of the derived def\/ining relations and the set of the def\/ining ordering relation have a controllable ``limiting behavior'' is introduced in Subsection~\ref{seclimit}. Subsection~\ref{subsection_sl_n} deals with the diagonal reduction algebra for $\sl_n$; the quadratic Casimir operator for
${\mathrm{DR}}(\sl_n)$ as well as for the diagonal reduction algebra for an arbitrary semi-simple
Lie algebra~$\k$ is given there. Subsection~\ref{subsection3.4} is devoted to the stabilization and cut phenomena with respect to the embedding of the Lie algebra $\gl_n\oplus\gl_m$ into
the Lie algebra~$\gl_{n+m}$; the theorem about the behavior of the centers of the diagonal reduction algebra under the cutting is proved there.

\subsection{New variables}

We shall use the following elements of $\Uh$:
\begin{gather*}  \tphi_{ij}:=\frac{\th_{ij}}{\th_{ij}-1}  ,\quad \tpsi_{ij}:=\frac{\th_{ij}-1}{\th_{ij}} ,
\quad \ttphi_{ij}:=\frac{\th_{ij}-1}{\th_{ij}-2} ,\quad \ttpsi_{ij}:=\frac{\th_{ij}-2}{\th_{ij}-1} ,\quad \Cprime_{ij}:=\frac{\th_{ij}-3}{\th_{ij}-2},
\end{gather*}
the variables $\th_{ij}$ are def\/ined in~\rf{hrond}. Note that $\tphi_{ij}\tpsi_{ij}=\ttphi_{ij}\ttpsi_{ij}=1$.

Def\/ine  elements $\tt_1,\ldots, \tt_n\in\Z_n$ by
\begin{gather*}
\label{3.1}
\tt_1:=t_1 ,\quad\tt_2:=\q_{1}(t_1) ,\quad\tt_3:=\q_{2}\q_{1}(t_1) ,
\quad\ldots, \quad \tt_n:=\q_{n-1}\cdots \q_{2}\q_{1}(t_1)  .
\end{gather*}
Using \rf{not7} we f\/ind the relations
 \begin{gather}
 \begin{split}
& \q_i(t_i)  =-\frac{1}{\th_{i,i+1}-1}t_i+\frac{\th_{i,i+1}}{\th_{i,i+1}-1}t_{i+1}  ,
\qquad \q_i(t_{i+1})= \frac{\th_{i,i+1}}{\th_{i,i+1}-1}t_i-\frac{1}{\th_{i,i+1}-1}t_{i+1} ,\\ 
& \q_i(t_k) =  t_k ,\qquad k\not=i,i+1  ,
 \end{split}\label{acwot}
 \end{gather}
which can be used to convert the def\/inition \rf{3.1} into  a linear over the ring $\Uh$ change of variables:
\begin{gather}
\begin{split}
& \ds{\tt_l} =\ds{t_l\prod_{j=1}^{l-1}\tphi_{jl}
-\sum_{k=1}^{l-1}t_k\frac{1}{\th_{kl}-1}\prod_{j=1}^{k-1}\tphi_{jl}}  , \\
& \ds{t_l} =\ds{\tt_l\prod_{j=1}^{l-1}\tpsi_{jl} +\sum_{k=1}^{l-1}\tt_k\frac{1}{\th_{kl}}
\prod_{\substack{j=1 \\ j\neq k}}^{l-1}\tpsi_{jk}} .
\end{split} \label{3.1a}
\end{gather}
For example,
\begin{gather*}
\tt_2 =-\frac{1}{\th_{12}-1}t_1+\frac{\th_{12}}{\th_{12}-1}t_2 ,\qquad
t_2=\frac{1}{\th_{12}}\tt_1+\frac{\th_{12}-1}{\th_{12}}\tt_2 ,\\
\tt_3 =-\frac{1}{\th_{13}-1}t_1-\frac{\th_{13}}{(\th_{13}-1)(\th_{23}-1)}t_2+
\frac{\th_{13}\th_{23}}{(\th_{13}-1)(\th_{23}-1)}t_3  ,\\
t_3 =  \frac{\th_{12}+1}{\th_{12}\th_{13}}\tt_1+\frac{\th_{12}-1}{\th_{12}\th_{23}}\tt_2+
\frac{(\th_{13}-1)(\th_{23}-1)}{\th_{13}\th_{23}}\tt_3  .
\end{gather*}

In terms of the new variables $\tt$'s, the linear in $t$ central element (\ref{clit}) reads
\begin{gather*}
\sum t_i=\sum\tt_i\prod_{a:a\neq i}\frac{\th_{ia}+1}{\th_{ia}} .
\end{gather*}

\subsection{Braid group action}
\label{brgac}

Since
$\q_{i}^2(x)=x$ for any element $x$ of zero weight, the braid group acts as its symmetric group quotient on the space of weight 0 elements.
It follows from \rf{3.1} and $\q_{i}(t_1)=t_1$ for all $i>1$ that
\begin{gather}\label{3.2}
\q_{\sigma}(\tt_i)=\tt_{\sigma(i)}
\end{gather}
for any $\sigma\in \S_n$.

The action of the Zhelobenko automorphisms, see Section \ref{section2},  on the generators $\s_{kl}$  looks as follows:
\begin{alignat}{3}
& \q_i(\s_{ik})=-\s_{i+1,k}\tphi_{i,i+1} , \qquad  &&
\q_i(\s_{ki})=-\s_{k,i+1} ,\quad k\not=i,i+1 ,& \nonumber\\
&\q_i(\s_{i+1,k})=\s_{i,k} ,\qquad && \q_i(\s_{k,i+1})=\s_{k,i}\tphi_{i,i+1} ,\quad k\not=i,i+1 ,& \label{3.5}\\
& \q_i(\s_{i,i+1})=-\s_{i+1,i}\tphi_{i,i+1}\ttphi_{i,i+1} ,\qquad && \q_i(\s_{i+1,i})=-\s_{i,i+1} , & \nonumber\\
&\q_{i}(\s_{j,k})=\s_{j,k} ,\quad j,k\neq i,i+1 .\qquad &&&\nonumber
\end{alignat}

Denote $i'=n+1-i$, as before. The braid group action \rf{3.5} is compatible with the anti-involution $\epsilon$
and the involution $\omega$ (note that $\omega (\th_{ij})=\th_{j'i'}$),  see~\rf{anep} and~\rf{not2a}, in the following sense:
\begin{gather}
\epsilon  \q_i =\q_i^{-1}\epsilon  ,\\
 \omega  \q_i =\q_{i'-1}\omega .\label{qeps}
\end{gather}

Let $w_0$ be the longest element of the Weyl group of $\gl_n$, the symmetric group $\S_n$. Similarly to the squares of
the transformations corresponding to the simple roots, see \rf{invr}, the action of~$\q_{w_0}^2$ is the conjugation by a
certain element of $\Uh$.

\begin{lem}\label{lqw02}
We have
\begin{gather}\label{qw02}
\q_{w_0}^2(x)=S^{-1}xS  ,
\end{gather}
where
\begin{gather}
 S=\prod_{i,j:i<j}\th_{ij}  .\label{conjq02}
 \end{gather}
\end{lem}

The proof shows that the formula \rf{qw02} works for an arbitrary reductive Lie algebra, with $S=\prod_{\alpha\in \Delta_+}\th_{\alpha}$.

\begin{prop}\label{paqw0}
The action of $\q_{w_0}$ on generators reads
\begin{gather}
 \label{3.6}
 \q_{w_0}(\s_{ij})=(-1)^{i+j} \s_{i'j'}\prod_{a:a<i'}\tphi_{ai'}
\prod_{b:b>j'}\tphi_{j'b} ,\\
\q_{w_0}(\tt_i)=\tt_{i'} .\label{3.6a}
\end{gather}
\end{prop}

The proofs of Lemma \ref{lqw02} and Proposition \ref{paqw0} are in Section \ref{sectionproofs}.

\subsection{Def\/ining relations}\label{section3.3}

To save  space we omit in this section the symbol $\mult$ for the multiplication in the algebra~$\Z_n$. It should not lead to any confusion since no other multiplication is used in this section.

Each relation which we will derive will be of a certain weight, equal to a sum of two roots. From general considerations the upper estimate
for the number of terms in a quadratic relation of weight $\lambda=\alpha+\beta$ is the number $|\lambda|$ of quadratic combinations $\s_{\alpha'}\s_{\beta'}$
with $\alpha'+\beta'=\lambda$.
There are several types, excluding the trivial one, $\lambda=2(\ve_i-\ve_j)$, $|\lambda|=1$:
\begin{enumerate}\itemsep=0pt
\item $\lambda=\pm(2\ve_i-\ve_j-\ve_k)$, where $i,j$ and $k$ are pairwise distinct. Then $|\lambda|=2$.
\item $\lambda=\ve_i-\ve_j+\ve_k-\ve_l$ with pairwise distinct $i,j,k$ and $l$. Then $|\lambda|=4$.
\item $\lambda=\ve_i-\ve_j$, $i\neq j$. For $\s_{\alpha'}\s_{\beta'}$, there are $2(n-2)$ possibilities (subtype 3a) with $\alpha'=\ve_i-\ve_k$,
$\beta'=\ve_k-\ve_j$ or $\alpha'=\ve_k-\ve_j$, $\beta'=\ve_i-\ve_k$ with $k\neq i,j$ and $2n$ possibilities (subtype 3b) with $\alpha'=0$,
$\beta'=\ve_i-\ve_j$ or $\alpha'=\ve_i-\ve_j$, $\beta'=0$. Thus $|\lambda|=4(n-1)$.
\item $\lambda=0$. There are $n^2$ possibilities (subtype 4a) with $\alpha'=0$, $\beta'=0$ and $n(n-1)$ possibilities (subtype 4b) with
$\alpha'=\ve_i-\ve_j$,  $\beta'=\ve_j-\ve_i$, $i\neq j$. Here $|\lambda|=n(2n-1)$.
\end{enumerate}

Below we write down relations for each type (and subtype) separately. The relations of the types 1 and 2 have a simple form in terms of the original
generators $\s_{ij}$. To write the relations  of the types 3 and 4, it is convenient to renormalize the generators $\s_{ij}$
with $i\not=j$. Namely, we set
\begin{gather}
\label{not8}\hs_{ij}=\s_{ij}\prod_{k=1}^{i-1}\tphi_{ki}  .
\end{gather}

In terms of the generators $\hs_{ij}$, the formulas \rf{3.5} for the action of the automorphisms $\q_i$  translate as follows:
\begin{alignat*}{3}
& \q_{i}(\hs_{ik})=-\hs_{i+1,k} ,\qquad &&  \q_{i}(\hs_{i+1,k})=\hs_{i,k}\tphi_{i+1,i} ,\quad k\not=i,i+1 , & \\
&\q_{i}(\hs_{ki})=-\hs_{k,i+1} ,\qquad && \q_{i}(\hs_{k,i+1})=\hs_{k,i}\tphi_{i,i+1}=\tpsi_{i+1,i}\hs_{k,i} ,\quad k\not=i,i+1,&
\\
&\q_{i}(\hs_{i,i+1})=-\tpsi_{i+1,i}\hs_{i+1,i} , \qquad &&\q_{i}(\hs_{i+1,i})=-\hs_{i,i+1}\tphi_{i+1,i} ,& \\
 &\q_{i}(\hs_{j,k})=\hs_{j,k} , \quad j,k\neq i,i+1 .\qquad &&&
 \end{alignat*}

{\bf 1.} The relations of  the type 1 are:
\begin{gather}\label{relation1}
\s_{ij}\s_{ik}=\s_{ik}\s_{ij}\tphi_{kj} ,\qquad
\s_{ji}\s_{ki}=\s_{ki}\s_{ji}\tpsi_{kj} ,\qquad \text{for}\quad j<k ,\  i\not=j,k .
\end{gather}

{\bf 2.} Denote
\begin{gather*}
D_{ijkl}:=\left( {\ds\frac{1}{\th_{ik}}}-{\ds\frac{1}{\th_{jl}}}\right).
\end{gather*}
Then, for any four pairwise dif\/ferent indices $i$, $j$, $k$ and $l$, we have the following relations of the type 2:
\begin{gather}
\begin{split}
& [ \s_{ij},\s_{kl} ]=\s_{kj}\s_{il}D_{ijkl} , \qquad i<k ,\ j<l  , \\ 
& \s_{ij}\s_{kl}-\s_{kl}\s_{ij}\tphi_{jl}'\tphi_{lj}'=\s_{kj}\s_{il}
D_{ijkl} ,\qquad i<k ,\ j>l.
\end{split} \label{relation2}
\end{gather}

{\bf 3a.} Let $i\neq k\neq l\neq i$. Denote
\begin{gather*}
\mathring{E}_{ikl}:=-\left((\tt_i-\tt_k) \frac{\th_{il}+1}{\th_{ik}\th_{il}}+(\tt_k-\tt_l) \frac{\th_{il}-1}{\th_{kl}\th_{il}}\right)\hs_{il}+\sum_{a:a\neq i,k,l}
\hs_{al}\hs_{ia} \frac{\ttphi_{ai}}{\th_{ka}+1} .
\end{gather*}
With this notation the f\/irst group of the relations of the type 3 is:
\begin{gather}\notag\hs_{ik}\hs_{kl}\tpsi_{ik} -\hs_{kl}\hs_{ik}\ttphi_{ki}
=\mathring{E}_{ikl} , \qquad i<k<l ,\\   \notag
\hs_{ik}\hs_{kl}\tpsi_{ik}\tpsi_{lk}\ttphi_{lk} -\hs_{kl}\hs_{ik}\ttphi_{ki}
=\mathring{E}_{ikl}  , \qquad i<l<k ,\\
\label{relation3}
\hs_{ik}\hs_{kl}\tphi_{ki} -\hs_{kl}\hs_{ik}\ttphi_{ki}=\mathring{E}_{ikl} , \qquad  k<i<l ,\\
\notag \hs_{ik}\hs_{kl}\tphi_{ki}\tphi_{li}\ttpsi_{li} -\hs_{kl}\hs_{ik}\ttphi_{ki}
=\mathring{E}_{ikl}  ,\qquad  k<l<i ,\\
\notag \hs_{ik}\hs_{kl}\tpsi_{ik}\tpsi_{lk}\ttphi_{lk}\tphi_{li}\ttpsi_{li}-\hs_{kl}\hs_{ik}\ttphi_{ki}
 =\mathring{E}_{ikl}  , \qquad l<i<k , \\  \notag
\hs_{ik}\hs_{kl}\tphi_{ki}\tpsi_{lk}\ttphi_{lk}\tphi_{li}\ttpsi_{li}
-\hs_{kl}\hs_{ik}\ttphi_{ki}=\mathring{E}_{ikl} , \qquad l<k<i .
\end{gather}

The relations \rf{relation3} can be written in a more compact way with the help of both systems, $\s_{ij}$ and $\hs_{ij}$, of generators. Let now
\begin{gather*}
E_{ikl}:=-\left((\tt_i-\tt_k) \frac{\th_{il}+1}{\th_{ik}\th_{il}}+(\tt_k-\tt_l) \frac{\th_{il}-1}{\th_{kl}\th_{il}}\right)
\s_{il}+\sum_{a:a\neq i,k,l}
\hs_{al}\s_{ia} \frac{\ttphi_{ai}}{\th_{ka}+1} .
\end{gather*}
Then
\begin{gather}
\begin{split}
& \s_{ik}\hs_{kl}\tpsi_{ik} -\hs_{kl}\s_{ik}\ttphi_{ki}=E_{ikl} ,\qquad k<l,\\ 
& \s_{ik}\hs_{kl}\tpsi_{ik}\tpsi_{lk}\ttphi_{lk} -\hs_{kl}\s_{ik}\ttphi_{ki}=E_{ikl} ,\qquad l<k .
\end{split}
\label{shfo3a}
\end{gather}
Moreover, after an extra redef\/inition: $\!\mathring{\hspace{.1cm}z}_{kl}\!\!\!\!\!\!\mathring{\phantom{z}}\ =\hs_{kl}\ttphi_{lk}$ for $k>l$, the left hand side of the
second line in \rf{shfo3a} becomes, up to a common factor, the same as the left hand side of the f\/irst line, namely, it reads $(\s_{ik}
\mathring{\hspace{.1cm}z}_{kl}\!\!\!\!\!\!\mathring{\phantom{z}}\ \ \tpsi_{ik} -\mathring{\hspace{.1cm}z}_{kl}\!\!\!\!\!\!\mathring{\phantom{z}}\ \ \s_{ik}\ttphi_{ki})
\tpsi_{lk}$.

{\bf 3b.} Let $i\neq j\neq k\neq i$. The second group of relations of the type 3 reads:
\begin{gather}\notag\hs_{ij}\tt_{i}=  \tt_i\hs_{ij}\Cprime_{ji}-\tt_j\hs_{ij}  \frac{1}{\th_{ij}+2}-
\sum_{a:a\not=i,j}\hs_{aj}\hs_{ia}  \frac{1}{\th_{ia}+2}  ,\\
\label{3.12}\hs_{ij}\tt_{j}= -\tt_i\hs_{ij}  \frac{\Cprime_{ji}}{\th_{ij}-1}+
\tt_j\hs_{ij}\tphi_{ij}\tpsi_{ji}\ttphi_{ji}+
\sum_{a:a\not=i,j}\hs_{aj}\hs_{ia}\tphi_{ij}\tpsi_{ji}{\ttphi_{ai}}{\th_{ja}+1}  ,\\  \notag
\hs_{ij}\tt_{k}= \tt_i\hs_{ij}  \frac{(\th_{ij}+3)\ttphi_{ji}}{(\th_{ik}^2-1)(\th_{jk}-1)}+
\tt_j\hs_{ij}  \frac{(\th_{ij}+1)\ttphi_{ji}}{(\th_{ik}-1)(\th_{jk}-1)^2}+
\tt_k\hs_{ij}\tphi_{ik}\tphi_{ki}\tphi_{jk}\ttpsi_{jk}\\
  \notag
  \phantom{\hs_{ij}\tt_{k}=}{}
- \hs_{kj}\hs_{ik}  \frac{(\th_{ij}+1)\ttphi_{ki}}{(\th_{ik}-1)(\th_{jk}-1)}
-\sum_{a:a\not=i,j,k}\hs_{aj}\hs_{ia}  \frac{\th_{ij}+1}{(\th_{ik}-1)(\th_{jk}-1)}
\frac{\ttphi_{ai}}{\th_{ka}+1}  .
\end{gather}

{\bf 4a.} The relations of the weight zero (the type~4) are also divided into 2 groups.
This is the f\/irst group of the relations:
\begin{gather}
 \label{relation4a}[\tt_i,\tt_j]=0.
 \end{gather}
As follows from the proof, the relations \rf{relation4a} hold for the diagonal reduction algebra for an arbitrary reductive
Lie algebra: the images of the generators, corresponding to the Cartan subalgebra, commute.

{\bf 4b.} Finally, the second group of the relations of the type 4 is
\begin{gather}
 \label{relation4}
 [\hs_{ij},\hs_{ji}]=\th_{ij}-\frac{1}{\th_{ij}}(\tt_i-\tt_j)^2+
\sum_{a:a\not=i,j}\left(\frac{1}{\th_{ja}+1}\hs_{ai}\hs_{ia}-\frac{1}{\th_{ia}+1}\hs_{aj}\hs_{ja}\right) ,
\end{gather}
where $i\neq j$.

{\bf Main statement.} Denote by $\mathfrak{R}$  the system \rf{relation1},  \rf{relation2}, \rf{relation3},
\rf{3.12},  \rf{relation4a} and \rf{relation4} of the relations.

\begin{theorem} \label{mthe}
The relations $\mathfrak{R}$ are the defining relations for the weight generators $\s_{ij}$ and
$t_i$ of the algebra $\Z_n$. In particular, the set  \eqref{not6} of ordering relations follows over $\Uh$ from $($and is equivalent to$)$  $\mathfrak{R}$.
\end{theorem}

The derivation of the system $\mathfrak{R}$ of the relations is given in Section~\ref{sectionproofs}. The validity in $\Z_n$ of relations from the
set~$\mathfrak{R}$, together with the results
from~\cite{KO3}, completes the proof of Theorem~\ref{mthe} (Section~\ref{proofmain}).

\subsection{Limit}\label{seclimit}

Let $\mathfrak{R}^\prec$ be the set of ordering relations~\rf{not6}.
Denote by ${\mathfrak R}_0$ the homogeneous (quadratic) part of the system ${\mathfrak R}$ and by $\mathfrak{R}^\prec_0$ the
homogeneous part of the system $\mathfrak{R}^\prec$.

{\bf 1.} Placing coef\/f\/icients from $\Uh$ in all relations from $\mathfrak{R}_0$ to the same side (to the right, for example) from the monomials
$\pp_L\mult\pp_M$, one can give arbitrary numerical values to the variables~$h_\alpha$ ($\alpha$'s~are roots of~$\f$).

The structure of the extremal  projector $P$ or the recurrence relation \rf{proof3}  implies that the system ${\mathfrak R}_0$ admits, for an arbitrary
reductive Lie algebra, the limit at $h_{\alpha_i}=c_i h$, $h\rightarrow\infty$ ($\alpha_i$ ranges through the set of simple positive
roots of  $\f$ and $c_i$ are generic positive constants). Moreover, this homogeneous algebra becomes the usual commutative (polynomial) algebra in
this limit; so this limiting behavior of the system $\mathfrak{R}_0$, used in the proof, generalizes to a wider class of reduction algebras, related to a pair $(\g,\k)$ as in the introduction.

{\bf 2.} The limiting procedure from paragraph {\bf 1} establishes the bijection between the set of relations and the set of
unordered pairs $(L,M)$, where $L,M$ are indices of basic vectors of $\p$.
 The proof in \cite{KO3} shows that over  $\UUh$ the system $\mathfrak{R}$ can be rewritten
in the form of ordering relations for an arbitrary order on the set $\{\pp_L\}$ of generators. Here $\UUh$ is  the f\/ield of fractions of the ring $\U (\h )$.

 By def\/inition, the relations from $\mathfrak{R}^\prec$ are labeled by pairs $(L,M)$
with $L>M$. The above bijection induces therefore a bijection between the sets $\mathfrak{R}$ and $\mathfrak{R}^\prec$.

\subsection[$\sl_n$]{$\boldsymbol{\sl_n}$}\label{subsection_sl_n}

{\bf 1.} Denote the subalgebra of $\Z_n$, generated by two central elements~\rf{clih} and~\r{clit}, by~$Y_n$; the algebra~$Y_n$ is isomorphic to~$\Z_1$.

Since the extremal projector for $\sl_n$ is the same as for $\gl_n$, the diagonal reduction algebra ${\mathrm{DR}}(\sl_n)$ for $\sl_n$ is
naturally a subalgebra of~$\Z_n$. The subalgebra ${\mathrm{DR}}(\sl_n)$ is complementary to~$Y_n$ in the sense
that $\Z_n=Y_n\otimes {\mathrm{DR}}(\sl_n)$.

The algebra ${\mathrm{DR}}(\sl_n)$ is generated by $\s_{ij}$, $i,j=1,\dots ,n$, $i\neq j$, and $t_{i,i+1}:=t_i-t_{i+1}$, $i=1,\dots ,n-1$
(and the Cartan subalgebra $\h$, generated by~$h_{i,i+1}$, of the diagonally embedded~$\sl_n$). The elements~$t_{i,i+1}$ form a basis
in the space of ``traceless'' combinations~$\sum c_mt_m$ (traceless means that  $\sum c_m=0$), $c_m\in \Uh$.

{\bf 2.} The action of the braid group restricts onto the traceless subspace:
\begin{alignat*}{3}
& \q_i(t_{i-1,i})=t_{i-1,i}+
 \frac{\th_{i,i+1}}{\th_{i,i+1}-1} t_{i,i+1} ,\qquad &&
\q_i(t_{i+1,i+2})= \frac{\th_{i,i+1}}{\th_{i,i+1}-1} t_{i,i+1}+t_{i+1,i+2} ,& \\
& \q_i(t_{i,i+1})=- \frac{\th_{i,i+1}+1}{\th_{i,i+1}-1}t_{i,i+1} ,\qquad & &
\q_i(t_{k,k+1})=t_{k,k+1} ,\quad  k\neq i-1,i,i+1 . &
\end{alignat*}

The traceless subspace with respect to the generators $t_i$ and the traceless subspace with respect to the generators $\tt_i$
(that is, the space of linear combinations $\sum c_m\tt_m$, $c_m\in \Uh$, with $\sum c_m=0$) coincide.  Indeed, in the expression of $t_l$ as a linear
combination of $\tt_k$'s (the second line in~\rf{3.1a}), we f\/ind, calculating residues and the value at inf\/inity, that the sum of the coef\/f\/icients is~1,
\begin{gather*}
\prod_{j=1}^{l-1}\tpsi_{jl}+\sum_{k=1}^{l-1} \frac{1}{\th_{kl}}
\prod_{\substack{j=1 \\ j\neq k}}^{l-1}\tpsi_{jk}=1  .
\end{gather*}
Therefore, in the decomposition of the dif\/ference $t_i-t_j$ as a linear combination of $\tt_k$'s, the sum of the coef\/f\/icients vanishes, so it is traceless with
respect to $\tt_k$'s; $t_{l,l+1}$ is a linear combination of $\tt_{12},\tt_{23},\dots ,\tt_{l,l+1}$ (and vice versa). It should be however noted that in contrast to~\rf{3.1a}, the coef\/f\/icients in these combinations do not factorize into a product of linear monomials, the lowest example
is $\tt_{34}$:
\begin{gather*}
\tt_{12}= \frac{\th_{12}}{\th_{12}-1}t_{12} ,\qquad \tt_{23}=\frac{\th_{23}}{\th_{13}-1}\left( -\frac{1}{\th_{12}-1}t_{12}
+\frac{\th_{13}}{\th_{23}-1}t_{23}\right) ,\\
\tt_{34}= \frac{\th_{34}}{\th_{14}-1}\left( -\frac{1}{\th_{13}-1}t_{12}-
\frac{\th_{14}(\th_{13}-1)+\th_{23}(\th_{24}-1)}
{(\th_{13}-1)(\th_{23}-1)(\th_{24}-1)}t_{23}
+\frac{\th_{14}\th_{24}}{(\th_{24}-1)(\th_{34}-1)}t_{34}\right)  .
\end{gather*}

{\bf 3.} One can directly see that the commutations between $\s_{ij}$ and the dif\/ferences $t_k-t_l$ close. The renormalization~\rf{not8} is compatible with
the $\sl$-condition and, as we have seen, the set $\{ t_{i,i+1}\}$ of generators can be replaced by the set $\{ \tt_{i,i+1}\}$. Therefore, one can work with the
generators $\hs_{ij}$, $i,j=1,\dots ,n$, $i\neq j$, and $\tt_{i,i+1}:=t_i-t_{i+1}$, $i=1,\dots ,n-1$. A direct look at the relations \rf{relation1},  \rf{relation2},
\rf{relation3},  \rf{3.12},  \rf{relation4a} and \rf{relation4} shows that the only non-trivial verif\/ication concerns the relations  \rf{3.12}; one has to check here
the following assertion: when $\hs$ moves through $\tt_{i,i+1}$, only traceless combinations of $\tt_l$'s appear in the right hand side. Write a~relation
from the list  \rf{3.12} in the form $\hs_{ij}\tt_l=\sum_m \chi_m^{(i,j,l,m)}\tt_m\hs_{ij}+\cdots$, $\chi_m^{(i,j,l,m)}\in \Uh$, where dots stand for terms with
$\hs\hs$. The assertion follows from the direct observation that for all~$i$,~$j$ and~$l$ the sum of the coef\/f\/icients $\chi_m^{(i,j,l,m)}$ is~1, $\sum_m \chi_m^{(i,j,l,m)}=1$.

{\bf 4.} With the help of the central elements \rf{clih}, \rf{clit} and \rf{drclit} one can build a unique linear in $t$'s traceless combination:
\begin{gather*}
\sum_{i=1}^n (h_i-2i)t_i-\left( \frac{1}{n} \sum_{i=1}^n h_i-n-1\right)
\sum_{j=1}^n t_j  .
\end{gather*}
It clearly depends only on the dif\/ferences $h_i-h_j$ and belongs therefore to the center of the subalgebra ${\mathrm{DR}}(\sl_n)$.

One can write this central element in the form
\begin{gather}\label{eq4.24}
\sum_{u,v=1}^{n-1}{\mathsf{C}}^{uv}h_{u,u+1}t_{v,v+1}+\sum_{v=1}^{n-1}(n-v)vt_{v,v+1}
=\sum_{u,v=1}^{n-1}{\mathsf{C}}^{uv}(\th_{u,u+1}+1)t_{v,v+1}  ,
\end{gather}
where ${\mathsf{C}}^{uv}$ is the inverse Cartan matrix of~$\sl_n$.

In general, let $\k$ be a semi-simple Lie algebra of rank $r$ with the Cartan matrix $a_{ij}$. Let $b_{ij}$ be the symmetrized Cartan matrix and $(\ ,\, )$
the  scalar product on $\h^*$ induced by the invariant non-degenerate bilinear form on $\k$, so that
\[
a_{ij}=d_i b_{ij} ,\qquad b_{ij}=(\alpha_i,\alpha_j) ,\qquad d_i =2/(\alpha_i,\alpha_i) .
\]
For each $i=1,\dots , r$ let $\alpha_i^\vee$ be the coroot vector corresponding to the simple root $\alpha_i$, so that $\alpha_j(\alpha_i^\vee)=a_{ij}$. Let $d_{ij}$ be the matrix, inverse to $c_{ij}=d_i b_{ij}d_j
$. Let $\rho\in\h^*$ be the half-sum of all positive roots. Write
\[
\rho=\frac{1}{2}\sum_{i=1}^r  n_i\alpha_i  ,
\]
where $n_i$ are nonnegative integers. Let $t_{\alpha_i}$ be the images of $H_{\a_i}={\alpha_i^\vee}^{(1)}-{\alpha_i^\vee}^{(2)}$ in the diagonal
reduction algebra  ${\mathrm{DR}}(\k)$ and $h_{\a_i}={\alpha_i^\vee}^{(1)}+{\alpha_i^\vee}^{(2)}$ be the coroot vectors of the diagonally
embedded Lie algebra $\k$. The generalization of the central element~\eqref{eq4.24} to the reduction algebra ${\mathrm{DR}}(\k)$ reads
\[
\sum_{i,i=1}^r d_{ij} h_{\alpha_i} t_{\alpha_j}+\sum_{i=1}^r n_i (\alpha_i,\alpha_i) t_{\alpha_i} .
\]

\subsection{Stabilization and cutting} \label{subsection3.4}

In \cite{KO3} we discovered the stabilization and cut phenomena which are heavily used in our derivation of the set of def\/ining relations for the
diagonal reduction algebras of $\gl$-type. The consideration in \cite{KO3} uses the standard (by the f\/irst coordinates) embedding of $\gl_n$ into $\gl_{n+1}$.
In this subsection we shall make several more precise statements about the stabilization and cut considering now the embedding of $\gl_n\oplus\gl_1$
into $\gl_{n+1}$ (more generally, $\gl_n\oplus\gl_m$ into $\gl_{n+m}$). These precisions are needed to establish the behavior of the center of
the diagonal reduction algebra: namely we shall see that cutting preserves the centrality.

Notation: $\h$ in this subsection denotes the Cartan subalgebra of $\gl_{n+m}$.

Consider an embedding of $\gl_n\oplus \gl_m$ into $\gl_{n+m}$, given by an assignment $e_{ij}\mapsto e_{ij}$, $i,j=1,\ldots,n$, and
$e_{ab}\mapsto e_{n+a,n+b}$, $a,b=1,\ldots,m$, where $e_{kl}$ in the source are the generators of $\gl_n\oplus \gl_m$ and target $e_{kl}$
are in $\gl_{n+m}$. This rule together with the similar rule $E_{ij}\mapsto E_{ij}$ and $E_{ab}\mapsto E_{n+a,n+b}$ def\/ines an embedding of the Lie
algebra $(\gl_n\oplus \gl_m)\oplus(\gl_n\oplus \gl_m)$ into the Lie algebra $\gl_{n+m}\oplus\gl_{n+m}$ and of the enveloping algebras
$\Ar_n\otimes \Ar_m=\U(\gl_n\oplus\gl_n)\otimes \U(\gl_m\oplus\gl_m)$ into $\Ar_{n+m}=\U(\gl_{n+m}\oplus\gl_{n+m})$. This embedding clearly
maps nilpotent subalgebras of $\gl_n\oplus \gl_m$ to the corresponding nilpotent subalgebras of $\gl_{n+m}$ and thus def\/ines an embedding
$\iota_{n,m}:\Z_n\otimes\Z_m\to \Z_{n+m}$ of the corresponding double coset spaces. However, the map $\iota_{n,m}$ is not a~homomorphism of
algebras. This is because the multiplication maps are def\/ined with the help of projectors, which are dif\/ferent for $\gl_n\oplus \gl_m$ and $\gl_{n+m}$.

However, as we will explain now we can control certain dif\/ferences between the two multiplication maps. Let $\V_{n,m}$ be the left ideal of the
algebra $\Z_{n+m}$ generated by elements $\s_{ia}$ with $i=1,\ldots,n$ and $a=n+1,\ldots, n+m$; let  $\V_{n,m}'$ be the right ideal of the
algebra $\Z_{n+m}$ generated by elements $\s_{ai}$ with $i=1,\ldots, n$ and $a=n+1,\ldots, n+m$.

Write any element $\lambda\in Q_+$ (the positive cone of the root lattice of $\gl_{n+m}$) in the form $\lambda=\sum_{k=1}^{n+m}\lambda_k\varepsilon_k$.
The element $\lambda$ can be presented as a sum
\begin{gather}\label{lalaplapp}\lambda= \lambda'+\lambda''  ,\end{gather}
where $\lambda'$ is an element of the root lattice of $\gl_n\oplus \gl_m$, and $\lambda''$ is
proportional to the simple root $\varepsilon_n-\varepsilon_{n+1}$: $\lambda'=\sum_{k=1}^{n+m}\lambda'_k\varepsilon_k$ with
$\sum_{k=1}^n\lambda'_k =  \sum_{k=n+1}^{n+m}\lambda'_k=0$ and $\lambda''=c(\varepsilon_n-\varepsilon_{n+1})$.

\begin{lem}\label{lemma1}
The left ideal  $\V_{n,m}\subset\Z_{n+m}$ consists of images in $\Z_{n+m}$ of  sums
$\sum_{ia}X_{ia} E_{ia}$ with $X_{ia}\in\Ab_{n+m}$, $i=1\lcd n$ and $a=n+1,\dots, n+m$.

The right ideal   $\V_{n,m}\subset\Z_{n+m}$ consists of images in $\Z_{n+m}$ of  sums $\sum_{ai} E_{ai}Y_{ai}$ with $Y_{ai}\in\Ab_{n+m}$,
$i=1\lcd n $ and  $a=n+1
\lcd n+m$.\end{lem}

\begin{proof}
Present the projector $P$ for the Lie algebra $\gl_{n+m}$ as a sum of terms
\begin{gather}
\label{wedeopro} \xi e_{-\gamma_1}\cdots e_{-\gamma_t}e_{\gamma'_{1}}\cdots e_{\gamma'_{t'}}  ,
\end{gather}
where $\xi\in\Uh$, $\gamma_1,\dots ,\gamma_t$ and $\gamma'_1,\dots ,\gamma'_{t'}$ are positive roots of $\gl_{n+m}$. For any $\lambda\in Q_+$
 denote by $P_\lambda$  the sum of above elements with $\gamma_1+\dots +\gamma_{t}=
\gamma'_1+\dots +\gamma'_{t'}=\lambda$. Then $P=\sum_{\lambda\in Q_+}P_\lambda$.
{} For any $X,Y\in\Ab$ def\/ine the element  $X\mult_\lambda Y$ as the image of $XP_\lambda Y$ in the reduction algebra. We have
$X\mult Y=\sum_{\lambda\in Q_+} X\mult_\lambda Y$.

For any $X\in\Ab_{n+m}$, $i=1\lcd n$ and $a=n+1\lcd n+m$ consider the product $X\mult_\lambda \s_{ia}$.

The product $X\mult_\lambda \s_{in}$ is zero if $\lambda''\not=0$ (the component $\lambda''$ is def\/ined by~(\ref{lalaplapp})).
Indeed, in this case in each summand of $P_\lambda$ one of $e_{\gamma'_{k'}}$ is equal to some $e_{jb}$, $j=1\lcd n$ and $b= n+1\lcd n+m$.
Choose an ordered basis of $\n_+$ which ends by all such $e_{jb}$ (ordered arbitrarily);
any element of~$\U(\n_+)$ can be written as a sum of ordered monomials, that is, monomials in which all such~$e_{jb}$ stand on the right.
Since $[e_{jb},E_{ia}]=0$ for any $i,j=1\lcd n$ and $a,b=n+1\lcd n+m$, the product
$e_{\gamma'_{k'}}E_{ia}$ belongs to the left ideal $\Ib$ and thus $X\mult_\lambda \s_{ia}=0$ in $\Z_{n+m}$.

If $\lambda''=0$ then generators of $\n_+$ in monomials entering the decomposition of  $P_\lambda$ are among the elements $e_{ij}$,
$1\leq i<j\leq n$, and $e_{ab}$, $n+1\leq a<b\leq n+m$ and thus their adjoint action leaves the space, spanned by all $E_{ia}$, $i=1\lcd n$,
$a=n+1\lcd n+m$ invariant, so $X\mult_\lambda \s_{ia}$ can be  presented as an image of the sum $\sum_{jb}X_{jb} E_{jb}$ with $X_{jb}\in\Ab_{n+m}$,
$j=1\lcd n$, $b=n+1\lcd n+m$. Thus, the left ideal, generated by all $\s_{ia}$ is contained in the vector space of images in $\Z_{n+m}$ of  sums
$\sum_{jb}X_{jb} E_{jb}$.

Moreover, for any $X\in\Ab_{n+m}$ the element $X\mult \s_{ia}$ is the image of $XE_{ia}+\sum_{j,b:\, j<i,\, b>a}\, X^{(jb)}E_{jb}$ for some $X^{(jb)}$
and the double induction on $i$ and $a$ proves the inverse inclusion.

The second part of lemma is proved similarly.
\end{proof}

\begin{cor}\label{corollary1}
We have the following decomposition of the free left $($and right$)$ $\Uh$-modules:
\begin{gather}\label{st4}
\Z_{n+m}= \mathrm{I}_{n,m}\oplus \Uh\cdot\iota_{n,m}(\Z_n\otimes\Z_m)   ,
\end{gather}
where $\mathrm{I}_{n,m}:=\V_{n,m}+\V'_{n,m}$.
\end{cor}

\begin{proof} The double coset space $\Z_{n+m}$ is a free left and right $\Uh$-module with a basis consisting of images of ordered monomials on elements
$E_{ij}$, $i,j=1\lcd n+m$; recall that we always use orders compatible with the partial order $<$ on $\h^*$, see (c) in Section~\ref{section2},
paragraph~2. We can choose an order for which all ordered monomials are of the form $XYZ$, where $X$ is a monomial on $E_{ai}$ with
$i=1\lcd n$ and $a=n+1\lcd n+m$, $Z$ is a monomial on $E_{ia}$  with $i=1\lcd n$ and $a=n+1\lcd n+m$ while $Y$ is a monomial on
$E_{ij}$ with $i,j=1\lcd n$ or $i,j=n+1\lcd n+m$. Then we apply the lemma above.
\end{proof}

For a moment denote for each $k>0$ the multiplication map in $\Z_k$ by $\mmult_{(k)}:\Z_k\otimes\Z_k\to\Z_k$  (instead of the default notation $\mult$,
see~\rf{not5a}); denote also for each $k,l>0$ by $\mmult_{(k,l)}$ the multiplication map $\mult_{(k)}\otimes \mult_{(l)}$ in $\Z_k\otimes \Z_l$.

\begin{prop}
\label{proposition1}
For any $x,y\in  \Z_n\otimes\Z_m$ we have
\begin{gather*}
\iota_{n,m}(x)\mmult_{(n+m)}\iota_{n,m}(y) = \iota_{n,m}(x\mmult_{(n,m)}y)+ z  ,
\end{gather*}
where $z$ is some element of $\mathrm{J}_{n,m}:=\V_{n,m}\cap\V'_{n,m}$.
\end{prop}

Let $\h_n$ and $\h_m$ be the Cartan subalgebras of $\gl_n$ and $\gl_m$, respectively. Denote the space
$\Z_n\otimes_{\overline{\mathrm{U}}(\h_n)}\Uh\otimes_{\overline{\mathrm{U}}(\h_m)}\Z_m$ by $\Uh\cdot\left(\Z_n\otimes\Z_m\right)$. The composition law
$\mmult_{(n,m)}$ naturally extends to the space $\Uh\cdot\left(\Z_n\otimes\Z_m\right)$ equipping it with an associative algebra structure (we keep the
same symbol $\mmult_{(n,m)}$ for the extended composition law in $\Uh\cdot\left(\Z_n\otimes\Z_m\right)$). Also, the map $\iota_{n,m}$ admits a natural
extension to a map $\iota_{n,m} : \Uh\cdot\left(\Z_n\otimes\Z_m\right)\to\Z_{n+m}$ denoted by the same symbol and def\/ined by the rule
$\iota_{n,m}(\varphi x):=\varphi\,\iota_{n,m}(x)$ for any $\varphi\in\Uh$ and $x\in \Z_n\otimes\Z_m$. The statement of Proposition~\ref{proposition1} remains
valid for this extension as well, that is, one can take $x,y\in\Uh\cdot\left(\Z_n\otimes\Z_m\right)$ in the formulation.

\begin{proof}[Proof  of Proposition~\ref{proposition1}.] Denote by $P_{n,m}:=P_n\otimes P_m$ the projector for the Lie algebra $\gl_{n}\oplus \gl_m$.

It is suf\/f\/icient to prove the following statement. Suppose $X$ and $Y$ are (non-commutative) polynomials in $E_{ij}$ with $i,j=1\lcd n$
Then the product of $x$ and $y$ in $\Z_{n+m}$ coincides with the image in $\Z_{n+m}$ of $X\, P_{n,m} Y$ modulo the left ideal $\V_{n,m}$
and modulo the right ideal $\V'_{n,m}$).

Due to the structure of the projector the condition $\lambda''=0$, see (\ref{lalaplapp}), implies that
the product $X\mult_\lambda Y$ related to $\gl_n\oplus \gl_m$ coincides with product $X\mult_\lambda Y$ related to $\gl_{n+m}$.

Let now $\lambda''\not=0$. Then
each monomial $e_{\gamma'_{1}}\cdots e_{\gamma'_{t'}}$ in the decomposition of $P_\lambda$, see~(\ref{wedeopro}), contains generators $e_{ia}$ with $i\in\{1\lcd n\}$ and  $a\in\{ n+1\lcd n+m\}$; these $e_{ia}$ can be assumed to be right factors of the corresponding monomial
(like in the proof of Lemma~\ref{lemma1}). The commutator of any such generator $e_{ia}$ with
every factor in $Y$ is a linear combination of the elements $E_{jb}$ with $j\in\{ 1\lcd n\}$ and
$b\in\{ n+1\lcd n+m\}$. Moving the resulting $E_{jb}$ to the right we see that the product $X\mult_\lambda Y$ is the image in $\Z_{n+m}$ of an element
of the form $\sum_s\, X_sY_s$ where each $Y_s$ belongs to the left ideal of $\Ab_{n+m}$ generated by $E_{jb}$ with
$j\in\{ 1\lcd n\}$ and $b\in\{ n+1\lcd n+m\}$ (one can say more: each $Y_s$ can be written in a form $\sum_{j,b}  Y_s^{(jb)}
E_{jb}$ where each $Y_s^{(jb)}\in\Ab_{n+m}$ does not involve generators $E_{ck}$ with $k\in\{ 1\lcd n\}$ and $c\in\{ n+1\lcd n+m\}$; we don't need
this stronger form). Thus, due to Lemma~\ref{lemma1}, $X\mult_\lambda Y\in\V_{n,m}$.

Similarly, each $X_s$ participating in the sum $\sum_s  X_sY_s$, see above, belongs to the right ideal of~$\Ab_{n+m}$ generated by the elements $E_{bj}$ with
$j\in\{ 1\lcd n\}$ and $b\in\{ n+1\lcd n+m\}$. So, again by Lemma~\ref{lemma1}, $X\mult_\lambda Y\in\V'_{n,m}$.
\end{proof}

Suppose that we have a relation
\begin{gather}
\sum_k a_k\mmult_{(n,m)} b_k=0 , \label{st1}
\end{gather}
where all $a_k$ and $b_k$ are elements of $Z_n\otimes Z_m$. Then, due to Proposition~\ref{proposition1}, we have the following relation in $\Z_{n+m}$:
 \begin{gather}
 \sum_k \bar{a}_k\mmult_{(m+n)} \bar{b}_k=z  , \label{st2}
\end{gather}
where $\bar{a}_k=\iota_{n,m}(a_k)$, $\bar{b}_k=\iota_{n,m}(b_k)$ and $z\in \mathrm{J}_{n,m}=\V_{n,m}\cap\V'_{n,m}$.

On the other hand, suppose we have the following relation in $\Z_{n+m}$:
\begin{gather}
 \sum_k \bar{a}_k\mmult_{(m+n)} \bar{b}_k=u  , \label{st3}
\end{gather}
where all $a_k$ and $b_k$ are elements of $Z_n\otimes Z_m$, $\bar{a}_k=\iota_{n,m}(a_k)$, $\bar{b}_k=\iota_{n,m}(b_k)$, and
$u\in \mathrm{I}_{n,m}=\V_{n,m}
+\V'_{n,m}$. Then the elements $a_k$ and $b_k$ satisfy the relation \rf{st1} and  $u\in \mathrm{J}_{n,m}$. Indeed, suppose that the relation \rf{st3} is
satisf\/ied and $ \sum_k a_k\mmult_{(n,m)} b_k=v$ for some $v\in Z_n\otimes Z_m$. It follows from Proposition~\ref{proposition1} that $\sum_k \bar{a}_k\mmult_{(m+n)} \bar{b}_k-\bar{v}$ belongs to $\mathrm{J}_{n,m}$; here $\bar{v}=\iota_{n,m}(v)$.
Then~\rf{st3} implies that $\bar{v}\in \mathrm{I}_{n,m}$ and thus $\bar{v}=0$ due to Corollary~\ref{corollary1}. Thus $v=0$, since the map $\iota_{n,m}$ is an
inclusion, and $u\in\mathrm{J}_{n,m}$.

We refer to the implication \rf{st1} $\Rightarrow$ \rf{st2} as   {\it stabilization}.
Call {\it cutting} the (almost inverse) implication \rf{st3}~$\Rightarrow$~\rf{st1} which can be understood as a procedure of getting
relations in $\Z_{n}\otimes\Z_m$ from relations in $\Z_{n+m}$; we say that \rf{st1} is the {\it cut} of~\rf{st3}.
Clearly all relations in $\Z_{n}\otimes\Z_m$ can be obtained by cutting appropriate relations in~$\Z_{n+m}$.

Let $\pi_{n,m}:\Z_{n+m}\to\Uh\cdot\left(\Z_{n}\otimes\Z_m\right)$ be the composition of the projection $\bar{\pi}_{n,m}$ of $\Z_{n+m}$ onto
$\iota_{n,m}(\Uh\cdot\Z_n\otimes\Z_m)=\Uh\cdot\iota_{n,m}(\Z_n\otimes\Z_m)$ along $\mathrm{I}_{n,m}$, see~\rf{st4},
and of the inverse to the inclusion~$\iota_{n,m}$:
\[
\pi_{n,m}=\iota_{n,m}^{-1}\circ\bar{\pi}_{n,m}  .
\]
We have the following consequence of Proposition~\ref{proposition1} and Corollary~\ref{corollary1}.

\begin{prop}
\label{proposition2}
Let $x$ be a central element of $\Z_{n+m}$. Then $\pi_{n,m}(x)$ is a central element of
$\Uh\cdot\left(\Z_{n}\otimes\Z_m\right)$.
\end{prop}

\begin{proof}
Denote $X=\pi_{n,m}(x)$. Then, by def\/inition, $x=\iota_{n,m}(X)+z$, where $z\in\mathrm{I}_{n,m}$. Since $x$ is central, it is of zero weight;
so $X$ and $z$ are of zero weight as well. Thus each monomial entering the decomposition of $z$ contains both types of generators, $E_{ai}$ and
$E_{ia}$, where $i\in\{1\lcd n\}$ and $a\in\{ n+1\lcd n+m\}$, which implies that $z\in \mathrm{J}_{n,m}=\V_{n,m}\cap\V'_{n,m}$.
Take any $Y\in\Z_{n}\otimes\Z_m$. We now prove that $X\mmult_{(n,m)}Y-Y\mmult_{(n,m)}X=0$. Denote $y=\iota_{n,m}(Y)$.
Due to Proposition \ref{proposition1},
\begin{gather}\label{st5}
\iota_{n,m}(X\mmult_{(n,m)}Y-Y\mmult_{(n,m)}X)=(x-z)\mmult_{(m+n)}y-y\mmult_{(m+n)}(x-z)+z'  ,
\end{gather}
where $z'\in \mathrm{J}_{n,m}=\V_{n,m}\cap\V'_{n,m}$. Since $x$ is central in $\Z_{n+m}$, the right hand side of~\rf{st5}
is equal~to
\[
y\mmult_{(m+n)}z-z\mmult_{(m+n)}y+z'  ,
\]
which is an element of $\mathrm{I}_{n,m}=\V_{n,m}\oplus \V'_{n,m}$ since $z,z'\in \mathrm{J}_{n,m}$. On the other hand, the left hand side of~\rf{st5}
belongs to $\Uh\cdot\iota_{n,m}(\Z_n\otimes\Z_m)$. Thus, by Corollary~\ref{corollary1}, both sides of~\rf{st5} are equal to zero and
$X\mmult_{(n,m)}Y-Y\mmult_{(n,m)}X=0$ since the map $\iota_{n,m}$ is injective.
\end{proof}

The map $\pi_{n,m}$ obeys properties similar to those of the Harish-Chandra map $\U(\g)^{\h}\to \U(\h)$ ($\U(\g)^{\h}$ is the space of elements of zero
weight). For instance, its restriction to the center of~$\Z_{n+m}$ is a homomorphism. More precisely, if $x$ is a central element of $\Z_{n+m}$, then
\begin{gather}
\begin{split}
&   \pi_{n,m}(x\mmult_{(m+n)}y)=\pi_{n,m}(x) \mmult_{(n,m)}\pi_{n,m}(y)\qquad \text{and}\\ 
& \pi_{n,m}(y\mmult_{(m+n)}x)=\pi_{n,m}(y) \mmult_{(n,m)}\pi_{n,m}(x)
\end{split} \label{st5bis}
\end{gather}
for any $y\in\Z_{n+m}$. Indeed, let $X=\pi_{n,m}(x)$, $Y=\pi_{n,m}(y)$. Then
\[
x=\iota_{n,m}(X)-z  ,\qquad y=\iota_{n,m}(Y)-u  ,
\]
where $u\in  \mathrm{I}_{n,m}$ while, as it was noted in the proof of Proposition~\ref{proposition2}, $z\in  \mathrm{J}_{n,m}$.
Moreover, it is clear that $z$ can be written in the form $z=\sum_az'_az^{\phantom{'}}_a$, where $z_a\in\V_{n,m}$ and $z'_a\in\V'_{n,m}$ (for instance,
use the order as in the proof of Corollary~\ref{corollary1}). Then (dropping for brevity the multiplication symbol $\mmult_{(m+n)}$) we have
\begin{gather}
 \iota_{n,m}(X)\iota_{n,m}(Y)= (x+z)(y+u)
 =\left(x+\sum_az'_az^{\phantom{'}}_a\right)(y+\tilde{z}'+\tilde{z})\nonumber\\
\phantom{\iota_{n,m}(X)\iota_{n,m}(Y)}{}
 =xy+\sum_a z'_az^{\phantom{'}}_a(y+\tilde{z}'+
\tilde{z})+x\tilde{z}+\tilde{z}'x  \equiv  xy\ \mod \mathrm{I}_{n,m}  .\label{hcpfce}
\end{gather}
Here $\tilde{z}\in\V_{n,m}$ and $\tilde{z}' \in\V'_{n,m}$. In the last equality we used the centrality of~$x$. Due to Proposi\-tion~\ref{proposition1},~(\ref{hcpfce}) is
precisely equivalent to the f\/ist part of~\rf{st5bis}. The second part of~\rf{st5bis} is proved similarly.

\section{Proofs}\label{sectionproofs}

\subsection[Tensor $\J$]{Tensor $\boldsymbol{\J}$}\label{sectionJ}

The multiplication map $\mult$ in $\Z_n$ (we return to the original notation) is given by the
prescription~\rf{not5a}, as in any reduction algebra. It can be formally expanded into
a series over the root lattice of certain bilinear maps as follows. Set
\begin{gather*}
\UU(\bb_\pm):=\Uh\otimes_{\U(\h)}\U(\bb_\pm)  , \qquad \UUb:=\UU(\bb_-)\otimes_{\Uh}  \UU(\bb_+)  .
\end{gather*}
All these are associative algebras. Besides, both algebras~$\UU(\bb_\pm)$ are $\Uh$-bimodules. The
algebra~$\UUb$ admits three commuting actions of~$\Uh$. Two of them are given by the assignments
\begin{gather*}
X(Y\otimes Z):=XY\otimes Z  ,\qquad (Y\otimes Z)X:=Y\otimes ZX  ,
\end{gather*}
for any $X\in\Uh$, $Y\in\UU(\bb_-)$ and $Z\in\UU(\bb_+)$.
The third action associates to any $X\in\Uh$, $Y\otimes Z\in\UUb$ the element $YX\otimes Z=Y\otimes XZ\in\UUb$.

Present the projector $P$ in an ordered form:
\begin{gather} \label{proof1}
P=
\sum_{\gamma,i}\grave{F}_{\gamma,i}\grave{E}_{\gamma,i}\grave{H}_{\gamma,i}=
\sum_{\gamma,i}\grave{H}_{\gamma,i}\grave{F}_{\gamma,i}\grave{E}_{\gamma,i}  ,
\end{gather}
the summation is over $\gamma\in\Q_+$ and $i\in {\mathbb{Z}}_{\geq 0}$;
every $\grave{F}_{\gamma,i}$ is an element of $\U(\n_-)$ of the weight~$-\gamma$, every
$\grave{E}_{\gamma,i}$ is an element of $\U(\n_+)$ of the weight $\gamma$ and
$\grave{H}_{\gamma,i}\in\Uh$. Let $\J$ be the following element of $\UUb$:
\begin{gather*}
\J :=\sum_{\gamma,i}\grave{F}_{\gamma,i}\otimes\grave{E}_{\gamma,i}\grave{H}_{\gamma,i}=
\sum_{\gamma,i}\grave{H}_{\gamma,i}\grave{F}_{\gamma,i}\otimes\grave{E}_{\gamma,i}  ,\qquad \gamma\in\Q_+  ,\quad i\in {\mathbb{Z}}_{\geq 0}  .
\end{gather*}
Due to the PBW theorem in $\U(\gl_n)$ the tensor~$\J$ is uniquely def\/ined by the projector $P$; it is
of total weight zero: $h\J=\J h$ for any $h\in\h$. We have the weight decomposition of $\J$ with
respect to the adjoint action of $\h$ in the second tensor factor of~$\UUb$:
\begin{gather*}
\J= \bigoplus_{\lambda\in\Q_+}\J_\lambda  ,
\end{gather*}
where $\J_\lambda$ consists of all the terms, corresponding to $\grave{F}_{\lambda,i}\grave{E}_{\lambda,i}
\grave{H}_{\lambda,i}$ in~\rf{proof1}
(contributing to $\lambda\in\Q_+$ in the summation),
\begin{gather*}
\J_\lambda :=\sum_i\grave{F}_{\lambda,i}\otimes\grave{E}_{\lambda,i}\grave{H}_{\lambda,i}  .
\end{gather*}
By def\/inition of $\J$, the multiplication $\mult$ in the double coset space $\Z_n$ can be described by the relation
\begin{gather}\label{proof4}
a\mult b= m\left( (a\otimes 1)\J (1\otimes b )\right),
\end{gather}
where $m(\sum_ic_i\otimes d_i)$ is the image in $\Z_n$ of the element $\sum_ic_id_i$.
Moreover, in \rf{proof4} we can replace all products $\grave{E}_{\gamma,i}b$ in the second tensor
factor by the adjoint action of $\grave{E}_{\gamma,i}$ on $b$ (in fact, for $\grave{E}_{\gamma,i}=e_{\gamma_m}\cdots e_{\gamma_1}$,
we can replace $\grave{E}_{\gamma,i}b$ by $[\grave{E}_{\gamma,i},b]$ or by $\hat{e}_{\gamma_m}\cdots \hat{e}_{\gamma_1}(b)$,
see \rf{dehaco}) and likewise all products $a\grave{F}_{\gamma,i}$
in the f\/irst tensor factor by the opposite adjoint action of $\grave{F}_{\gamma,i}$ on $a$.
We have a~decomposition of the product $\mult$ into a sum over $\Q_+$:
\begin{gather}\label{proof5}
a\mult b=\sum_{\lambda\in\Q_+}a\mult_\lambda b  ,\qquad\text{where}\quad
a\mult_\lambda b:=m\left( (a\otimes 1)  \J_\lambda (1\otimes b )\right)  .
\end{gather}
If $a$ and $b$ are weight elements of $\Z_n$ of weights $\nu(a)$ and $\nu(b)$, then the product
$a\mult_\lambda b$ is the image in $\Z_n$ of the sum $\sum_i a_ib_i$, where the weight of each
$b_i$ is $\nu(b)+\lambda$, and the weight of each~$a_i$ is $\nu(a)-\lambda$.

The tensor $\J$ satisf\/ies the Arnaudon--Buf\/fenoir--Ragoucy--Roche (ABRR) dif\/ference equation~\cite{ABRR}, see also~\cite{K} for the translation of the results of~\cite{ABRR} to the language
of reduction algebras. To describe the equation, let  $\vartheta=\frac{1}{2}
\sum_{k=1}^n\th_k^2\in \U(\h)$;
for any positive root $\gamma\in\Delta_+$, denote by~$T_\gamma$ the following linear operator on the vector space~$\UUb$:
\begin{gather*}T_{\gamma}(X\otimes Y):=X e_{-\gamma}\otimes e_{\gamma}Y  .\end{gather*}
The ABRR equation means the relation \cite{ABRR,K}:
\begin{gather*}
[1\otimes \vartheta, \J]=-\sum_{\gamma\in\Delta_+}T_{\gamma}(\J)
 .
\end{gather*}
This relation is equivalent to the following system of recurrence relations for the weight components
$\J_\lambda$:
\begin{gather} \label{proof3}
\J_\lambda\cdot\left(\th_\lambda+\frac{(\lambda,\lambda)}{2}\right)=
-\sum_{\gamma\in\Delta_+}T_\gamma\left(\J_{\lambda-\gamma}\right)  ,
\end{gather}
where $\th_\lambda :=\sum_k\lambda_k\th_k$ for $\lambda=\sum_k\lambda_k\ve_k$. The recurrence
relations~(\ref{proof3})  together with the initial condition $\J_0=1\otimes 1$
uniquely determine all weight components~$\J_\lambda$.

It should be noted that the recurrence relations \rf{proof3} provides less information about the structure of the denominators (from~$\U (\h )$)
of the summands of the extremal projector $P$ than the information implied
 by the product formula (see~\cite{AST}) for the extremal projector.

Using \rf{proof3} we get in particular:
\begin{gather} \label{proof3a}
\J_{\a}  = -(\th_{\a}+1)^{-1}e_{-a}\otimes e_{\a} ,
\qquad \a=\ve_i-\ve_{i+1} ,
\\ \notag
 \J_{\a+\beta}  =  (\th_{\a+\b}+1)^{-1}\Bigl(-e_{-\a-\b}\otimes e_{\a+\b}+
(\th_{\a}+1)^{-1} e_{-\a}e_{-\b}\otimes e_{\b}e_{\a}\\
\phantom{\J_{\a+\beta}  =}{}  +(\th_{\b}+1)^{-1} e_{-\b}e_{-\a}\otimes e_{\a}e_{\b}\Bigr)
  ,\qquad \a=\ve_{i-1}-\ve_{i}  ,\quad \b=\ve_{i}-\ve_{i+1}  ,\label{proof3b}
\\
  \label{proof3c}\J_{\ve_i-\ve_{j}+\ve_k-\ve_{l}}= \J_{\ve_i-\ve_{j}}\cdot
\J_{\ve_k-\ve_{l}}  ,\qquad  i<j<k<l  .
\end{gather}

\subsection{Braid group action}

The proof of the relations \rf{acwot} and \rf{3.5} consists of the following arguments, valid for any
reduction algebra. Let $\a$ be any simple root of $\gl_n$, $\a=\ve_i-\ve_{i+1}$ and $\g_\a$ the
corresponding $\sl_2$ subalgebra of $\gl_n$. It is spanned by the elements $e_\a=e_{i,i+1}$,
$e_{-\a}= e_{i+1,i}$ and $h_\a= h_i-h_{i+1}$. Let $\sy_\a=\sy_i$ be the corresponding
automorphism of the algebra $\Ar$ and $\q_\a=\q_i$ the Zhelobenko automorphism of
$\Z_n$. Assume that $Y\in\Ar$ belongs, with respect to the adjoint action of $\g_\a$, to an irreducible
f\/inite-dimensional $\g_\a$-module of dimension $2j+1$, $j\in\{\ts0,1/2,1,\ldots\ts\}$. Assume further
that $Y$ is homogeneous, of weight $2\ts m$, $[h_\a,Y]=2m Y$. Identify $Y$ with its image in $\Z_n$. Then
$\q_\a(Y)$ coincides with the image in $\Z_n$ of the element
\begin{gather*}
\prod_{i=1}^{j+m}\,(h_\a+i+1)\ts\cdot\ts\sy_\alpha(Y)\cdot
{\ds\prod_{i=1}^{j+m} (h_\a-i+1)^{-1}}  .
\end{gather*}
This can be checked directly using \cite[Proposition 6.5]{KO}.

In the realization of irreducible $\sl_2$-modules as the spaces of homogeneous polynomials in two
variables $u$ and $v$,
\begin{gather*}e_\a\mapsto u\frac{\partial}{\partial v}  , \qquad h_\a\mapsto u\frac{\partial}{\partial u}-v\frac{\partial}{\partial v}
\qquad {\mathrm{and}}\qquad  e_{-\a}\mapsto v\frac{\partial}{\partial u}  ,
\end{gather*}
the operator $\sy_\alpha$ becomes $(\sy_\alpha f)(u,v)=f(-v,u)$, or, in the basis $|j,k\rangle :=x^{j+k}y^{j-k}$ ($j$ labels the
representation; $k=0,1,\dots ,2j$),
\begin{gather*}\sy_\alpha : \  |j,-j+k\rangle\mapsto (-1)^k |j,j-k\rangle  .\end{gather*}

\begin{proof}[Proof  of Lemma \ref{lqw02}, Subsection~\ref{brgac}.] 
To see this, write a reduced expression for $\q_{w_0}$, $\q_{w_0}=\q_{\alpha_{i_1}}\cdots \q_{\alpha_{i_M}}$
with $\alpha_{i_1},\dots ,\alpha_{i_M}$ simple roots. Then $\q_{w_0}=\q_{\alpha_{i_M}}\cdots \q_{\alpha_{i_1}}$ as well.
Writing, for $\q_{w_0}^2$, the second expression after the f\/irst one, we get squares of
$\q_{\alpha_{i_s}}$'s (which are conjugations by $\th_{\alpha_{i_s}}^{-1}$'s; they thus commute) one after another.
Moving these conjugations to the left through the remaining~$\q$'s, we produce, exactly like in the
construction of a system of all positive roots from a reduced expression for
the longest element of the Weyl group of a reductive Lie group, the conjugation by the product (\ref{conjq02})
over all positive roots.
\end{proof}

\begin{proof}[Proof  of Proposition~\ref{paqw0}, Subsection~\ref{brgac}.] 
Only formula \rf{3.6} needs a proof (formula~\rf{3.6a}
is a particular case of~\rf{3.2}).

For a moment, denote the longest element
of the symmetric group $S_n$ by $\q_{w_0}^{(n)}$. Let $\psi_j:=\q_j\q_{j-1}\cdots\q_1$ (the product in the descending order).
We have $\q_{w_0}^{(n+1)}=\q_{w_0}^{(n)}\psi_n$ and  $\q_{w_0}^{(n+1)}=\psi_1\psi_2\cdots\psi_n$ (the product in the
ascending order).

For $j<n$ it follows from \rf{3.5} that $\psi_j (\s_{n+1,1})=(-1)^j\s_{n+1,j+1}$ (say, by induction on $j$). So,
\begin{gather*}\psi_n(\s_{n+1,1})=q_n\psi_{n-1} (\s_{n+1,j+1})=(-1)^{n-1}q_n(\s_{n+1,n})=(-1)^{n}\s_{n,n+1}  ,\end{gather*}
again by \rf{3.5}. Next, $\psi_k\psi_{k+1}\cdots\psi_{n-1}(\s_{n,n+1})=\s_{k,n+1}$ by induction on $n-k$ and
again \rf{3.5}. Thus,
\begin{gather}\label{binh}
\q_{w_0}(\s_{n+1,1})=(-1)^n \s_{1,n+1} ,
\end{gather}
establishing \rf{3.6} for $i=n+1$ and $j=1$. We now prove \rf{3.6} for $i>j$ (positions below the main diagonal) by induction backwards
on the height $i-j$ of a negative root;  the formula~\rf{binh} serves as the induction base. Assume that~\rf{3.6}
is verif\/ied for a given level $i-j$ and $i-j-1>0$ (so that the positions $(i,j+1)$ and $(i-1,j)$ are still under the main
diagonal). By \rf{3.5}, $\s_{i,j+1}=-\q_j(\s_{ij})$, therefore
\begin{gather*}\q_{w_0}(\s_{i,j+1}) = -\q_{w_0}(\q_j(\s_{ij}))=-\q_{j'-1}(\q_{w_0}(\s_{ij})) \\
\phantom{\q_{w_0}(\s_{i,j+1})}{}  = (-1)^{i+j+1} \q_{j'-1}\left(\s_{i'j'}\ds{\prod_{a:a<i'}}\tphi_{ai'}  \ds{\prod_{b:b>j'}}\tphi_{j'b}\right)\\
\phantom{\q_{w_0}(\s_{i,j+1})}{} = (-1)^{i+j+1}\s_{i',j'-1}\tphi_{j'-1,j'}\ds{\prod_{a:a<i'}}\tphi_{ai'}  \ds{\prod_{b:b>j'}}\tphi_{j'-1,b}\\
\phantom{\q_{w_0}(\s_{i,j+1})}{} = (-1)^{i+j+1}\s_{i',(j+1)'}\ds{\prod_{a:a<i'}}\tphi_{ai'}\ds{\prod_{b:b>(j+1)'}}\tphi_{(j+1)',b}  .
\end{gather*}
In the second equality we used the identity $\q_{w_0}\q_j=\q_{j'-1}\q_{w_0}$ in the braid group; the third equality is the
induction assumption; in the fourth equality we used that $i'\neq j'-1$ (since $i-j-1>0$) and then~\rf{3.5}; in the f\/ifth equality we
replaced $j'-1$ by $(j+1)'$. The calculation for $\q_{w_0}(\s_{i-1,j})$ is similar; it uses $\s_{i-1,j}=\q_{i-1}(\s_{ij})$.
The proof of the formula~\rf{3.6} for positions below the main diagonal is f\/inished.

The proof of  \rf{3.6} for $i<j$ (positions above the main diagonal) follows now from Lem\-ma~\ref{lqw02}. 
\end{proof}

\subsection{Derivation of relations}\label{sectionDerRel}

The set  of def\/ining relations in $\Z_n$ divides into several dif\/ferent types,
see Section~\ref{section3.3}. We prove the necessary amount of relations
of each type and get the rest by applying the transformations from the braid group as well as the anti-involution~$\epsilon$, see~(\ref{anep}).

We never use the automorphism $\omega$, def\/ined in \rf{not2a}, in the derivation of relations. However, the involution $\omega$ is compatible with our set
of relations in the sense explained in Section \ref{proofmain}.

In the following we denote by the symbol $ \equiv $ the equalities of elements from $\Ab$ modulo the sum $(\Jb+\Ib )$ of two ideals $\Jb$ and $\Ib$
def\/ined in the beginning of Section \ref{section2}. Moreover,
for any two elements $X$ and $Y$ of the algebra $\Ab$ we may regard the expressions $X\mult Y$ and $X\mult_\lambda Y$ as the sums
of elements from $\Ab$ def\/ined in \rf{proof4} and \rf{proof5}. The sum $X\mult_\lambda Y$  is f\/inite. By the construction,
all but a f\/inite number of terms in the product $X\mult Y$ belong to $(\Jb+\Ib )$. Unlike to the system of notation
adopted in Section~\ref{section2}, our proof of each relation in $\Z_n$ will use equalities in $\Ab$ taken modulo $(\Jb+\Ib )$.

We also use the notation $H_i$ for the element $E_{ii}\in\Ar$ and $H_{ij}=H_i-H_j=E_{ii}-E_{jj}$.

{\bf 1.} We f\/irst prove in $\Z_n$ the relation
\begin{gather}\label{proof7}
\s_{12}\mult\s_{13}=\s_{13}\mult\s_{12}  \frac{\th_{23}}{\th_{23}+1}  .
\end{gather}
Elements $\s_{12}$ and $\s_{13}$ are images in $\Z_n$ of $E_{12}$ and $E_{13}$. Consider the
product $E_{12}\mult_\lambda E_{13}$. Since the adjoint action of $\gl_n$ preserves the space
$\p$, see Section~\ref{section-notation}, this product is the sum of such monomials $E_{ij}E_{kl}$,
with coef\/f\/icients in $\Uh$, that (i): the weight $\ve_k-\ve_l$
of $E_{kl}$ is equal to the weight $\ve_{1}-\ve_3$ of $E_{13}$ plus $\lambda\in\Q_+$, while (ii): the
weight $\ve_i-\ve_j$ of $E_{ij}$ is equal to the weight $\ve_{1}-\ve_2$ of $E_{12}$ minus
$\lambda$. Assume that $E_{12}\mult_\lambda E_{13}\neq 0$. By (i), $\lambda =-\ve_1+\ve_3+\ve_k-\ve_l$
and it can be positive only if $k=1$ and $l\geq 3$. So, the condition (i) implies that either
$\lambda=0$ or $\lambda=\ve_3-\ve_l$ with $l>3$. The possibility $\lambda=\ve_3-\ve_l$, $l>3$, is excluded
by the condition (ii). Therefore, $\lambda=0$ and
\begin{gather}
\label{proof8}
E_{12}\mult E_{13}\equiv E_{12}E_{13} .
\end{gather}
Similarly, for $\lambda\in\Q_+$, which can non-trivially contribute to the product $E_{13}\mult E_{12}$,
the analogue of the condition (i) on the weight $\lambda$ gives the
restriction $\lambda =0$ or $\lambda=\ve_2-\ve_k$, $k>2$; the analogue of the condition (ii)
further restricts $\lambda$:  $\lambda=0$ or $\lambda=\ve_2-\ve_3$, so we have
\begin{gather*}
E_{13}\mult_{\ve_2-\ve_3} E_{12}\equiv -E_{13}e_{32}e_{23}  \frac{1}{\th_{23}+1}E_{12}\equiv
-E_{13}e_{32}e_{23}E_{12}  \frac{1}{\th_{23}}\equiv E_{12}E_{13}  \frac{1}{\th_{23}}  ,
\end{gather*}
since $\J_{\ve_2-\ve_3}=-e_{32}\otimes e_{23}(\th_{23}+1)^{-1}$ as it follows from the ABRR
equation, see~(\ref{proof3a}), or from the precise explicit expression for the projector~$P$,
see~\cite{AST}. Thus, since $E_{12}$ and $E_{13}$ commute in the universal enveloping algebra
\begin{gather}\label{proof9}E_{13}\mult E_{12}\equiv E_{13}E_{12}+E_{13}\mult_{\ve_2-\ve_3} E_{12}=
E_{12}E_{13}\left(1+\frac{1}{\th_{23}}\right),
\end{gather}
Comparing \rf{proof8} and \rf{proof9} we f\/ind \rf{proof7}.

Applying to \rf{proof7} the anti-involution $\epsilon$, see (\ref{anep}), we get the relation
\begin{gather}\label{proof10}
\s_{21}\mult\s_{31}=\s_{31}\mult\s_{21}  \frac{\th_{23}+1}{\th_{23}}  .
\end{gather}

The rest of the relations \rf{relation1} are obtained from~\rf{proof7} and~\rf{proof10}
by applying dif\/ferent transformations $\q_w$ from the Weyl group.

{\bf 2.} Now we prove in $\Z_n$ the relation
\begin{gather}\label{proof11}
\s_{13}\mult\s_{24}-\s_{24}\mult\s_{13}=
\left(\frac{1}{\th_{12}}-\frac{1}{\th_{34}}\right)\s_{23}\mult\s_{14}  .
\end{gather}
We begin by the proof of this relation in $\Z_4$. We proceed in the same manner as for the derivation of the
relation (\ref{proof7}),
\begin{gather*}E_{13}\mult E_{24}\equiv E_{13}E_{24}+E_{13}\mult_{\ve_1-\ve_2} E_{24} \\
\phantom{E_{13}\mult E_{24}}{} \equiv
E_{13}E_{24}-E_{13}e_{21}e_{12}  \frac{1}{\th_{12}+1} E_{24}
 \equiv E_{13}E_{24}+E_{23}E_{14} \frac{1}{\th_{12}}  ,\\
 E_{24}\mult E_{13}\equiv E_{24}E_{13}+E_{24}\mult_{\ve_3-\ve_4} E_{13}\\
 \phantom{E_{24}\mult E_{13}}{}  \equiv
E_{24}E_{13}-E_{24}e_{43}e_{34}  \frac{1}{\th_{34}+1} E_{13} \equiv E_{13}E_{24}+E_{23}E_{14}  \frac{1}{\th_{34}}  ,\\
 E_{23}\mult E_{14} \equiv E_{23}E_{14}  .
 \end{gather*}
Combining the three latter equalities we obtain \rf{proof11} in $\Z_4$.

The dif\/ference of the left and right hand sides of \rf{proof11} in $\Z_n$ is a linear combination of monomials in $\s_{ij}$ of
the total weight $\ve_1+\ve_2-\ve_3-\ve_4$. The weight is non-trivial, so the monomials can be only quadratic.
Due to the stabilization phenomenon, each monomial should contain $\s_{ij}$ with $i>4$ or $j>4$, but, by the weight arguments,
there is no such non-zero possibility, which completes the proof of the relation \rf{proof11} in $\Z_n$.

The rest of relations \rf{relation2} is then obtained by applications of the transformations from the braid group.

{\bf 3a.} We continue and derive in $\Z_4$ the relation (we remind the notation $t_{ij}:=\s_{ii}-\s_{jj}$,
see Section~\ref{section2},
and $H_{ij}=E_{ii}-E_{jj}$):
\begin{gather} \label{proof12}
\s_{23}\mult\s_{12}-\s_{12}\mult\s_{23}=t_{12}\mult\s_{13}
 \frac{1}{\th_{12}}+t_{23}\mult\s_{13}  \frac{1}{\th_{23}}-\s_{43}\mult\s_{14}  \frac{\th_{34}+1}{\th_{34}\th_{24}}  .
\end{gather}
Using \rf{proof3a}--\rf{proof3c}, we calculate, to obtain the result for $\Z_4$:
\begin{gather}E_{12}\mult E_{23} \equiv E_{12}E_{23}+ E_{12}\mult_{\ve_1-\ve_2} E_{23}+
E_{12}\mult_{\ve_1-\ve_2+\ve_3-\ve_4} E_{23} \nonumber\\
\phantom{E_{12}\mult E_{23}}{}
 \equiv E_{12}E_{23}-H_{12}E_{13} \frac{1}{\th_{12}}  ,\label{1223r1}\\
  E_{23}\mult E_{12} \equiv E_{23}E_{12}+ E_{23}\mult_{\ve_2-\ve_3} E_{12}+
E_{23}\mult_{\ve_2-\ve_4} E_{12}\notag\\
\phantom{E_{23}\mult E_{12}}{} \equiv E_{23}E_{12}+H_{23}E_{13}  \frac{1}{\th_{23}}-E_{43}E_{14} \frac{(\th_{23}-1)}{\th_{23}
\th_{24}} ,\label{1223r2}\\
  H_{12}\mult E_{13} \equiv H_{12} E_{13}+H_{12}\mult_{\ve_3-\ve_4} E_{13}\equiv
H_{12} E_{13}  ,\label{1223r3}\\
  H_{23}\mult E_{13} \equiv H_{23} E_{13}+H_{23}\mult_{\ve_3-\ve_4} E_{13}\equiv
H_{23} E_{13}+E_{43}E_{14}  \frac{1}{\th_{34}}  ,\label{1223r4}\\
  E_{43}\mult E_{14} \equiv E_{43}E_{14} .\label{1223r5}
  \end{gather}
Combining the above equalities and taking into account that $[E_{12},E_{23}]=e_{13}\equiv 0$,
we get~\rf{proof12} in~$\Z_4$. We could apply here the stability arguments (as we shall do in the sequel)
but we give some more details at this point to give a f\/lavor of how such derivations of relations work.
For the same, as \rf{1223r1}--\rf{1223r5}, calculations for $\Z_n$, the analogues of the conditions (i)
and (ii), see paragraph 1 of this subsection, restrict $\lambda$ to be of the form $\ve_1-\ve_2+\ve_3-\ve_k$, $k\geq 3$ for
\rf{1223r1}; $\ve_2-\ve_k$, $k\geq 2$ for \rf{1223r2}; $\ve_3-\ve_k$, $k\geq 3$ for~\rf{1223r3}
and \rf{1223r4}; $\ve_4-\ve_k$, $k\geq 4$ for~\rf{1223r5}. It follows, for, say, $n=5$, that the right
hand sides of \rf{1223r1}--\rf{1223r5} might be modif\/ied only by an addition of the term proportional to
$E_{53}E_{15}$; and this will be compensated by an addition of the term, proportional to
$\s_{53}\mult\s_{15}$ to the right hand side of~\rf{proof12}, since $E_{53}\mult E_{15}\equiv E_{53}E_{15}$;
the proportionality coef\/f\/icient is uniquely def\/ined.
This pattern clearly continues and we conclude that there is a relation in $\Z_n$ of the form
\begin{gather}\label{proof13}
\s_{23}\mult\s_{12}-\s_{12}\mult\s_{23}=t_{12}\mult\s_{13} \frac{1}{\th_{12}}+t_{23}\mult\s_{13}\frac{1}{\th_{23}}-
\sum_{k>3}\s_{k3}\mult\s_{1k}X_k  ,
\end{gather}
with certain, uniquely def\/ined, coef\/f\/icients $X_k\in\Uh$, $k=4\lcd n$, and already known $X_4=(\th_{34}+1)
{\th_{34}^{-1}\th_{24}^{-1}}$. To f\/ind $X_5,\dots ,X_n$, we apply to \rf{proof13} the
automorphisms $\q_k$, $k=4\lcd n-1$, which leave invariant the left hand side and the f\/irst two terms in the
right hand side of~\rf{proof13}. The uniqueness of the relation of the form~\rf{proof13}, together with the equality
$\q_k(\s_{k3}\mult\s_{1k})=\s_{k+1,3}\mult\s_{1,k+1}(\th_{k,k+1}+1)\th_{k,k+1}^{-1}$, imply the recurrence relation
$X_{k+1}=\q_k(X_k)\cdot(\th_{k,k+1}+1)\th_{k,k+1}^{-1}$ and we f\/ind
\begin{gather*}
X_k=\frac{1}{\th_{2k}}\prod_{j=3}^{k-1}\frac{\th_{jk}+1}{\th_{jk}}  .
\end{gather*}

After the renormalization~\rf{not8} and the change of variables~\rf{3.1a}, the derived relation becomes one of the
relations in the f\/irst line of~\rf{relation3}.

Applying the transformations from the braid group, we obtain the rest of the relations from the list~\rf{relation3}.

{\bf 3b.} We have the following equalities in $\Z_3$:
\begin{gather}\label{proof15}
\s_{12}\mult t_1=t_1\mult\s_{12}  {\ds\frac{\th_{12}+2}{\th_{12}+1}}-t_2\mult\s_{12}  {\ds\frac{1}{\th_{12}+1}}
-\s_{32}\mult\s_{13}  {\ds\frac{\th_{23}+1}{\th_{23}(\th_{13}+1)}}  ,\\
\label{proof16}
\s_{12}\mult t_2=-t_1\s_{12}  {\ds\frac{1}{\th_{12}+1}}+t_2\mult\s_{12}
 {\ds\frac{\th_{12}+2}{\th_{12}+1}}+\s_{32}\mult\s_{13}  {\ds\frac{\th_{13}+2}{\th_{23}(\th_{13}+1)}}  ,
\end{gather}
and the equality in $\Z_4$:
\begin{gather}\label{proof17}
\s_{12}\mult t_4=t_4\mult\s_{12}-\s_{42}\mult\s_{14}  {\ds\frac{\th_{12}+1}{(\th_{14}+1)\th_{24}}}  .
\end{gather}
Equalities \rf{proof15} and \rf{proof16} are the results of the following calculations for $\Z_3$,
using \rf{proof3a}--\rf{proof3c}, and of the commutativity $[H_1,E_{12}]=e_{12}\equiv 0$,
$[H_2,E_{12}]=-e_{12}\equiv 0$:
\begin{gather*}E_{12}\mult H_1 \equiv E_{12}H_1+H_{12}E_{12}  \frac{1}{\th_{12}+1}-
E_{32}E_{13}  \frac{\th_{12}}{(\th_{12}+1)(\th_{13}+1)}  , \\
 E_{12}\mult H_2 \equiv E_{12}H_2-H_{12}E_{12}  \frac{1}{\th_{12}+1}-
E_{32}E_{13}  \frac{1}{(\th_{12}+1)(\th_{13}+1)}, \\
 H_1\mult E_{12} \equiv H_1E_{12}  , \qquad
 H_2\mult E_{12} \equiv H_2E_{12}-E_{32}E_{13}  \frac{1}{\th_{23}} , \qquad
 E_{32}\mult E_{13} \equiv E_{32}E_{13}  .\end{gather*}
The derivation of \rf{proof17} can be done with the help of the following calculations for $\Z_4$:
\begin{gather}\label{pr14n}
E_{12}\mult H_4 \equiv E_{12}H_4+E_{42}E_{14}  \frac{1}{\th_{14}+1},\\
 H_4\mult E_{12} \equiv H_4E_{12}+E_{42}E_{14}  \frac{1}{\th_{24}}  ,\qquad
 E_{42}\mult E_{14} \equiv E_{42}E_{14}  .\nonumber
 \end{gather}
We shall make a comment about the line \rf{pr14n} only. Here one might expect, by  the analogues of the conditions (i) and (ii), see paragraph 1 of this subsection,
non-trivial contributions to $E_{12}\mult H_4$ from the weights 0, $\ve_1-\ve_2$, $\ve_1-\ve_3$ and $\ve_1-\ve_4$.
So we need, in addition to \rf{proof3a}--\rf{proof3c}, some information about
$\J_{\ve_1-\ve_4}$. It follows from the ABRR equation that $\J_{\ve_1-\ve_4}(\th_{14}+1)=-T_{\ve_1-\ve_2}(\J_{\ve_2-\ve_4})
-T_{\ve_1-\ve_3}(\J_{\ve_3-\ve_4})-T_{\ve_1-\ve_4}(\J_0)-T_{\ve_2-\ve_3}(\J_{\ve_1-\ve_2+\ve_3-\ve_4})-
T_{\ve_2-\ve_4}(\J_{\ve_1-\ve_2})-T_{\ve_3-\ve_4}(\J_{\ve_1-\ve_3})$. Since $e_{13}$ and $e_{12}$
commute with $H_4$, the parts $T_{\ve_2-\ve_4}(\J_{\ve_1-\ve_2})$ and $T_{\ve_3-\ve_4}(\J_{\ve_1-\ve_3})$
do not contribute; $e_{42}$ and $e_{43}$ commute with $E_{12}$, so the parts  $T_{\ve_1-\ve_2}(\J_{\ve_2-\ve_4})$
and $T_{\ve_1-\ve_3}(\J_{\ve_3-\ve_4})$ do not contribute either; $\J_{\ve_1-\ve_2+\ve_3-\ve_4}=
\J_{\ve_3-\ve_4}\J_{\ve_1-\ve_2}$ does not contribute again since $e_{12}$ commute with $H_4$.
Thus the only contribution is from $T_{\ve_1-\ve_4}(\J_0)$ and we quickly arrive at~\rf{pr14n}.

Applying the automorphism $\q_3$ of the algebra $\Z_4$ to the relation \rf{proof17}, see \rf{acwot} and~\rf{3.5}, we f\/ind
\begin{gather*}
\s_{12}\mult\big(t_3\th_{34}-t_4\big)=\big(t_3\th_{34}-t_4\big)\mult \s_{12}
-\s_{32}\mult\s_{13}  {\ds\frac{\th_{34}(\th_{12}+1)}{(\th_{13}+1)\th_{23}}}  .
\end{gather*}
We then add \rf{proof17} to this relation  and obtain the following relation in $\Z_4$:
\begin{gather}\label{proof18}
\s_{12}\mult t_3=t_3\mult\s_{12}
-\s_{32}\mult \s_{13}  {\ds\frac{(\th_{12}+1)}{(\th_{13}+1)\th_{23}}}-
\s_{42}\mult\s_{14}  {\ds\frac{\th_{12}+1}{(\th_{14}+1)\th_{24}\th_{34}}}  .
\end{gather}

The stabilization arguments for \rf{proof15}, \rf{proof16} and \rf{proof18} imply the existence of the following relations in $\Z_n$:
\begin{gather}\label{proof19}\s_{12}\mult t_1 =t_1\mult\s_{12}  {\ds\frac{\th_{12}+2}{\th_{12}+1}}-
t_2\mult\s_{12}  {\ds\frac{1}{\th_{12}+1}}+\sum_{k>2}\s_{k2}\mult\s_{1k}X_k^{(1)},\\
 \label{proof20}\s_{12}\mult t_2 =-t_1\mult\s_{12}  {\ds\frac{1}{\th_{12}+1}}+t_2\mult\s_{12}  {\ds\frac{\th_{12}+2}{\th_{12}+1}}
+\sum_{k>2}\s_{k2}\mult\s_{1k}X_k^{(2)},
\\
\label{proof21}\s_{12}\mult t_3 =t_3\mult\s_{12}-\s_{32}\mult\s_{13}  {\ds\frac{(\th_{12}+1)}{(\th_{13}+1)\th_{23}}}
+\sum_{k>3}\s_{k2}\mult\s_{1k}X_k^{(3)},
\end{gather}
where all $X_k^{(i)}$ belong to $\Uh$ and the initial $X_k^{(i)}$ are known:
\begin{gather*} X_3^{(1)} =-\frac{\th_{23}+1}{\th_{23}(\th_{13}+1)}  ,\qquad
X_3^{(2)} =\frac{\th_{13}+2}{\th_{23}(\th_{13}+1)}  ,\qquad
X_4^{(3)} =-\frac{\th_{12}+1}{(\th_{14}+1)\th_{24}\th_{34}}  .
\end{gather*}

By the braid group transformation laws, $X_{k+1}^{(i)}=\q_k\big(X_k^{(i)}\big)\cdot(\th_{k,k+1}+1)\th_{k,k+1}^{-1}$
with $k>2$ for $i=1,2$ and $k>3$ for $i=3$, so that
\begin{gather*}
X_k^{(1)}=-\frac{1}{\th_{1k}+1}\prod_{j=2}^{k-1}\frac{\th_{jk}+1}{\th_{jk}} ,\qquad
X_k^{(2)}=\frac{\th_{1k}+2}{(\th_{1k}+1)\th_{2k}}\prod_{j=3}^{k-1}
\frac{\th_{jk}+1}{\th_{jk}}  ,\\
X_k^{(3)}=-\frac{\th_{12}+1}{(\th_{1k}+1)\th_{2k}\th_{3k}}\prod_{j=4}^{k-1}
\frac{\th_{jk}+1}{\th_{jk}}  .
\end{gather*}

After the renormalization \rf{not8} and the change of variables~\rf{3.1a}, the relations \rf{proof19}--\rf{proof21} turn into
the relations~\rf{3.12} for $i=1$, $j=2$ and $k=3$.

Applying the transformations from the braid group, we obtain the rest of the relations from the list \rf{3.12}.

{\bf 4a.} We now prove the relations \rf{relation4a} using the arguments similar to \cite[Subsection 6.1.2]{Zh}.
 Consider the
products $H_k\mult_\lambda H_l$  and $H_l\mult_\lambda H_k$ with $\lambda\not=0$.
These products are linear combinations, over $\Uh$, of monomials
\begin{gather*}
a_{kl;\vec{\gamma}}:=H_ke_{-\gamma_1}\cdots e_{-\gamma_m}e_{\gamma_m}\cdots e_{\gamma_1}H_l\qquad\text{and}\qquad
a_{lk;\vec{\gamma}}:=H_le_{-\gamma_1}\cdots e_{-\gamma_m}e_{\gamma_m}\cdots e_{\gamma_1}H_k  ,
\end{gather*}
respectively; here $m\geq 0$ and $\vec{\gamma}:=\{\gamma_1,\dots ,\gamma_m\}$. By construction, the coef\/f\/icient, from $\Uh$, of the monomial
$a_{kl;\vec{\gamma}}$ in $H_k\mult_\lambda H_l$ equals the coef\/f\/icient of $a_{lk;\vec{\gamma}}$ in $H_l\mult_\lambda H_k$.
The expressions $a_{kl;\vec{\gamma}}$ and $a_{lk;\vec{\gamma}}$ are both equal in $\Z_n$  to
\begin{gather*}(\gamma_1,\ve_k)(\gamma_1,\ve_l)E_{-\gamma_1}e_{-\gamma_2}\cdots
e_{-\gamma_m}e_{\gamma_m}\cdots e_{\gamma_2}E_{\gamma_1}  .\end{gather*}
Thus $H_k\mult_\lambda H_l\equiv H_l\mult_\lambda H_k$  for any $\lambda\not=0$. In the zero weight part $\mult_0$
of the product $\mult$ we have the equality $H_kH_l=H_lH_k$ as well. Therefore,
$H_k\mult H_l\equiv H_l\mult H_k$.

{\bf 4b.} The last group \rf{relation4} of relations is left. Like above, we f\/irst explicitly derive the following relation in~$\Z_3$:
\begin{gather} \s_{12}\mult \s_{21}=  h_{12}-t_{12}\mult t_{12}  \frac{1}{\th_{12}-1}+\s_{21}\mult\s_{12}
\frac{(\th_{12}-1)(\th_{12} +2)}{\th_{12}(\th_{12} +1)}\nonumber\\
\phantom{\s_{12}\mult \s_{21}=}{} +  \s_{31}\mult\s_{13} \frac{(\th_{12}-1) (\th_{13} +2)}{\th_{12}\th_{23}(\th_{13}+1)}
-\s_{32}\mult\s_{23}  \frac{\th_{23} +2}{(\th_{23} +1)\th_{13}}  \label{re01}
\end{gather}
(the f\/irst term in the right hand side is $h_{12}$, without hat). The relation \rf{re01} is a corollary of the following calculations
for $\Z_3$, together with the commutation relation $[E_{12},E_{21}]=h_{12}$,
\begin{gather}E_{12}\mult E_{21}\equiv   E_{12}E_{21}-H_{12}^2 \frac{1}{(\th_{12}-1)}
+E_{21}E_{12} \frac{2}{(\th_{12}-1)\th_{12}}\nonumber\\
\phantom{E_{12}\mult E_{21}\equiv}{} - E_{32}E_{23} \frac{\th_{12} -2}{(\th_{12}-1)\th_{13}}+
E_{31}E_{13} \frac{\th_{12} -2}{(\th_{12}-1)\th_{12}\th_{13}} ,\label{zemod1}\\
 H_{12}\mult H_{12}\equiv   H_{12}^2-E_{21}E_{12} \frac{4}{\th_{12} +1}
-E_{32}E_{23} \frac{1}{\th_{23} +1}\nonumber\\
\phantom{H_{12}\mult H_{12}\equiv}{} +  E_{31}E_{13}\left(-1+\frac{1}{\th_{23} +1}+\frac{4}{\th_{12} +1}\right)\frac{1}{\th_{13}+1} ,\label{zemod2}\\
E_{21}\mult E_{12}\equiv E_{21}E_{12}-E_{31}E_{13} \frac{1}{\th_{23}} ,\label{zemod3}\\
 E_{32}\mult E_{23}\equiv  E_{32}E_{23}-E_{31}E_{13} \frac{1}{\th_{12}} ,\label{zemod4}\\
 E_{31}\mult E_{13}\equiv  E_{31}E_{13} .\label{zemod5}
 \end{gather}
Here only the calculation of $E_{12}\mult E_{21}$ deserves a little explanation; by  the analogues of the conditions (i) and (ii),
see paragraph 1 of this subsection, non-trivial contributions to $E_{12}\mult E_{21}$ from the weights 0,
$\ve_1-\ve_2$, $2(\ve_1-\ve_2)$, $\ve_1-\ve_3$ and $2\ve_1-\ve_2-\ve_3$ are possible. By the ABRR equation,
$\J_{2(\ve_1-\ve_2)}(2\th_{12}+4)=-T_{\ve_1-\ve_2}(\J_{\ve_1-\ve_2})$ and $\J_{2\ve_1-\ve_2-\ve_3}
(2\th_1-\th_2-\th_3+3)=-T_{\ve_1-\ve_2}(\J_{\ve_1-\ve_3})-T_{\ve_1-\ve_3}(\J_{\ve_1-\ve_2})-T_{\ve_2-\ve_3}
(\J_{2(\ve_1-\ve_2)})$. We leave further details to the reader.

By the stabilization law in $\Z_4$ we have a relation
\begin{gather} \label{proof23}
\s_{12}\mult \s_{21}= h_{12}-t_{12}\mult t_{12}  \frac{1}{\th_{12}-1}
+\sum_{1\leq i<j\leq n}\s_{ji}\mult\s_{ij}X_{ij} ,\qquad X_{ij}\in\Uh
\end{gather}
with $n=4$, which dif\/fers from \rf{re01} by a presence of terms
\begin{gather*}\s_{43}\mult\s_{34} ,\qquad \s_{42}\mult\s_{24} ,\qquad \s_{41}\mult\s_{14} ,
\end{gather*}
with coef\/f\/icients in $\Uh$. Consider in $\Z_4$ the products  $\s_{12}\mult\s_{21}$, $t_{12}\mult t_{12}$
and $\s_{ji}\mult\s_{ij}$, $1\leq i<j\leq 4$. The weights $(\ve_3-\ve_4)-(\ve_i-\ve_j)$ do not belong
to the cone $\Q_+$ if $1\leq i<j<4$. Thus in the decomposition
\begin{gather*}
E_{ji}\mult E_{ij}\equiv \sum_{k<l}E_{lk}E_{kl}a_{kl} ,\qquad a_{kl}\in\Uh ,\quad 1\leq i<j<4 ,
\end{gather*}
the term with $E_{43}E_{34}$ has a zero coef\/f\/icient, $a_{34}=0$. The same statement holds for the products
$E_{41}\mult E_{14}$ and $E_{42}\mult E_{24}$ since the weights $(\ve_3-\ve_4)-(\ve_i-\ve_4)$ do
not belong to $\Q_+$ for $i=1,2$. Consider the product $E_{12}\mult E_{21}$. Here the term with
$E_{43}E_{34}$ is equal to $E_{12}\mult_{\ve_1-\ve_2+\ve_3-\ve_4}E_{21}$. By \rf{proof3c},
$\J_{\ve_1-\ve_2+\ve_3-\ve_4}=e_{43}e_{21}\otimes e_{12}e_{34}(\th_{12}+1)^{-1}(\th_{34}+1)^{-1}$ and
\begin{gather*}E_{12}\mult_{\ve_1-\ve_2+\ve_3-\ve_4}E_{21}=
E_{12}e_{43}e_{21}e_{12}e_{34}E_{21} \frac{1}{(\th_{12}-1)(\th_{34}+1)}\equiv 0,\end{gather*}
since $[e_{34},E_{21}]=0$ (and $[E_{12},e_{43}]=0$). In the similar manner, the term with
$E_{43}E_{34}$ in $H_{12}\mult H_{12}$ equals $H_{12}\mult_{\ve_3-\ve_4}H_{12}$ and vanishes
since $[e_{34},H_{12}]=0$.

On the other hand, the product $E_{43}\mult E_{34}$
def\/initely contains $E_{43}E_{34}=E_{43}\mult_0 E_{34}$. We thus conclude that the term
$\s_{43}\mult \s_{34}$ is absent in~\rf{proof23}, that is $X_{34}=0$.

For $n>4$, again by the stabilization law, we have a unique relation of the form \rf{proof23}. By uniqueness, it is invariant with respect to the
transformations $\q_3,\q_4,\dots ,\q_{n-1}$ which do not change the product $\s_{12}\mult \s_{21}$. Since $X_{34}=0$, we f\/ind, applying
$\q_4, \q_5, \dots ,\q_{n-1}$,  that  $X_{3j}=0$, $j>3$, wherefrom we further conclude, applying $\q_3,\q_4,\dots ,\q_{j-2}$, that
$X_{ij}=0$, $2<i<j$. We get f\/inally the following relation in $\Z_n$:
\begin{gather} \label{soo1221}
\s_{12}\mult \s_{21}=h_{12}-t_{12}\mult t_{12} \frac{1}{\th_{12}-1}
+\sum_{k=2\lcd n}\s_{k1}\mult\s_{1k}X_{1k}+\sum_{k=3\lcd n}\s_{k2}\mult\s_{2k}X_{2k}
\end{gather}
with known
\begin{gather*}
X_{12}=\frac{(\th_{12}-1)(\th_{12} +2)}{\th_{12}(\th_{12} +1)} ,\qquad
X_{13}=\frac{(\th_{12}-1) (\th_{13} +2)}{\th_{12}\th_{23}(\th_{13}+1)} ,\qquad
X_{23}=-\frac{\th_{23} +2}{(\th_{23} +1)\th_{13}} .
\end{gather*}
Applying to \rf{soo1221} the transformations $\q_3, \q_4,\dots ,\q_{n-1}$ we f\/ind by uniqueness
\begin{gather*}X_{i,k+1}=\frac{\th_{k,k+1}+1}{\th_{k,k+1}} \cdot\q_k(X_{ik}) ,\qquad i=1,2 ;\quad k=3,4,\dots ,n-1 ,
\end{gather*}
and thus
\begin{gather*}
 X_{1k}=\frac{(\th_{12}-1) (\th_{1k} +2)}{\th_{12}\th_{2k}(\th_{1k}+1)}\cdot\prod_{a=3}^{k-1}\frac{\th_{ak}+1}{\th_{ak}}  ,
\qquad X_{2k}=-\frac{\th_{2k} +2}{(\th_{2k} +1)\th_{1k}}\cdot\prod_{a=3}^{k-1}\frac{\th_{ak}+1}{\th_{ak}}
\end{gather*}
for $k>2$.

After the renormalization \rf{not8} and the change of variables~\rf{3.1a}, the relations~\rf{soo1221}  with the
obtained~$X_{1k}$ and $X_{2k}$ turns into the relation~\rf{relation4} for $i=1$ and $j=2$.

Applying the transformations from the braid group, we obtain the rest of the relations from the list~\rf{relation4}.

\subsection{Proof of Theorem \ref{mthe}}\label{proofmain}

For the proof of Theorem \ref{mthe} we just apply the results of~\cite{KO3}, which state that the system ${\mathfrak R}$ is the system of
def\/ining relations once it is satisf\/ied in the algebra $\Z_n$.

\begin{remark}
An attentive look shows that the system $\mathfrak{R}$ is closed under the anti-involution~$\epsilon$; that is, $\epsilon$~transforms any
relation from $\mathfrak{R}$ into a linear over $\Uh$ combination of relations from~$\mathfrak{R}$. Moreover, $\mathfrak{R}$ and
$\epsilon (\mathfrak{R})$ are equivalent over $\Uh$. Indeed, all relations in Section~\rf{sectionDerRel} were derived in
three steps: f\/irst we derive a relation in $\Z_n$ with $n\leq 4$; next by the stabilization principle we extend the derived relation
to~$\Z_n$ with arbitrary~$n$;
and then we f\/ind the whole list of relations of a given (sub)type by applying the braid group transformations (products of the
generators $\q_i$). Due to \rf{invr} one could use $\q_i^{-1}$ instead of $\q_i$ equivalently over
$\Uh$. A~straightforward calculation establishes the equivalence of the extended to arbitrary $n$ lists $\mathfrak{R}$ and
$\epsilon (\mathfrak{R})$ over $\Uh$ for $\Z_n$ with
$n\leq 4$ (this verif\/ication is lengthy for some relations).
Then with the help of~\rf{qeps} we f\/inish the check of the equivalence of $\mathfrak{R}$ and $\epsilon (\mathfrak{R})$ over $\Uh$ for $\Z_n$ with
arbitrary $n$.

Similar arguments establish  the equivalence of $\mathfrak{R}$ and $\omega (\mathfrak{R})$ over $\Uh$; here $\omega$ is the involution def\/ined
in~\rf{not2a}. In~\cite{KO3} this equivalence was obtained dif\/ferently, as a by-product of the equivalence, over $\Uh$, of the system $\mathfrak{R}$ and
the system  \rf{not6} of ordering relations.
\end{remark}

\section[Examples: $\sl_3$ and $\sl_2$]{Examples: $\boldsymbol{\sl_3}$ and $\boldsymbol{\sl_2}$}

In this section we write down the complete list of ordering relations for the diagonal reduction algebras
${\mathrm{DR}}(\sl_3)$ and ${\mathrm{DR}}(\sl_2)$. For
completeness we include the formulas for the action of the braid group generators and the expressions for the central elements.

We f\/irst give the list of relations for $\sl_3$. It is straightforward to give the list for $\sl_2$ directly; we comment however on how the list of relations and the
expressions for the central elements for $\sl(2)$ can be obtained by the cut procedure.

The list of relations for $\gl_3$ follows immediately from the list for $\sl_3$.

{\bf 1. Relations for $\boldsymbol{{\mathrm{DR}}(\sl_3)}$.} We write the ordering relations for the natural set of generators~$z_{ij}$, without redef\/initions.
We use here the following notation for $\sl_3$:
\begin{gather*}\z_\a :=\z_{12} ,\quad \z_\b :=\z_{23} ,\quad \z_{\a+\b} :=\z_{13} ,\quad \z_{-\a}:=\z_{21} ,
\quad \z_{-\b} :=\z_{23} ,\quad \z_{-\a-\b} :=\z_{31} ,
\\
t_\a:=t_{12} ,\quad t_\b:=t_{23} ,\quad h_\a:=h_{12} ,\quad h_\b:=h_{23} .
\end{gather*}
The relations are given for the following order $\grave{\succ}$ (this order was used in the proof in~\cite{KO3} of the completeness of
the set of relations):
\begin{gather}
 \z_{\a+\b}\ \grave{\succ}\ \z_{\a}\ \grave{\succ}\ \z_{\b}\ \grave{\succ}\ t_\b\ \grave{\succ}\ t_\a\ \grave{\succ}\ \z_{-\b}\ \grave{\succ}\ \z_{-\a}\
\grave{\succ}\ \z_{-\a-\b} .\label{stord}
\end{gather}

Due to the established in Theorem~\ref{mthe} and remarks in Section~\ref{seclimit} bijection between the set $\mathfrak{R}$ and the set $\mathfrak{R}^\prec$
of def\/ining relations, one can divide the ordering relations into the types, in the same way as we divided the def\/ining relations from the list~$\mathfrak{R}$.

The relations of type 1 are immediately rewritten as ordering relations:
\begin{gather}
\z_{\a +\b}\mult\z_\a =\z_\a\mult \z_{\a +\b}  \fracds{h_\b +2}{h_\b +1} ,\label{strel1}\\
\z_{\a +\b}\mult\z_\b =\z_\b\mult \z_{\a +\b} \fracds{h_\a +2}{h_\a +1} ,\\
\z_\a\mult \z_{-\b} =\z_{-\b}\mult\z_\a  \fracds{h_\a +h_\b +3}{h_\a +h_\b +2} ,\\
\z_\b\mult \z_{-\a} =\z_{-\a}\mult\z_\b  \fracds{h_\a +h_\b +3}{h_\a +h_\b +2} ,\\
\z_{-\a}\mult\z_{-\a -\b}=\z_{-\a -\b}\mult\z_{-\a} \fracds{h_\b +2}{h_\b +1} ,\\
\z_{-\b}\mult\z_{-\a -\b}=\z_{-\a -\b}\mult\z_{-\b} \fracds{h_\a +2}{h_\a +1} .
\end{gather}

The relations of type 2 are absent for $\sl_3$.

The ordering relations corresponding to the relations of type~3  we collect according to their weights. For each weight there
is one relation of subtype (3{\bf a}) and two relations of subtype~(3{\bf b}).

Weight $\alpha +\beta$:
\begin{gather}
\z_\a\mult\z_\b = -t_\a\mult \z_{\a +\b} \fracds{1}{h_\a +1}-t_\b\mult \z_{\a +\b} \fracds{1}{h_\b +1}+\z_\b\mult\z_\a  ,\\
 \nonumber \z_{\a +\b}\mult t_\a = t_\a\mult \z_{\a +\b}
 \fracds{h_\a h_\b+h_\b^2+2h_\a +6h_\b+9}{(h_\b +2)(h_\a +h_\b +3)}\\
\phantom{\z_{\a +\b}\mult t_\a =}{} +
t_\b\mult \z_{\a +\b} \fracds{h_\b^2+h_\a +6h_\b+9}{(h_\b +1)(h_\b +2)(h_\a +h_\b+3)}
 - \z_\b\mult \z_\a \fracds{h_\a +2h_\b +6}{(h_\a +2)(h_\b +2)} ,\\
 \nonumber \z_{\a +\b}\mult t_\b = t_\a\mult \z_{\a +\b} \fracds{h_\b}{(h_\b +2)(h_\a +h_\b+3)}+  t_\b\mult \z_{\a +\b}
 \fracds{h_\b \bigl( h_\a h_\b+h_\b^2+3h_\a +7h_\b+11\bigr)}{(h_\b +1)(h_\b +2)(h_\a +h_\b +3)}\\
\phantom{\z_{\a +\b}\mult t_\b =}{}
 +  \z_\b\mult \z_\a \fracds{2h_\a +h_\b +6}{(h_\a +2)(h_\b +2)} .
\end{gather}

Weight $\alpha$:
\begin{gather}
 \z_\a\mult t_\a =t_\a\mult \z_\a \fracds{h_\a +4}{h_\a +2}-\z_{-\b}\mult\z_{\a +\b} \fracds{h_\a +2h_\b +6}{(h_\b +1)(h_\a +h_\b +3)} ,\\
\z_\a\mult t_\b =-t_\a\mult \z_\a \fracds{1}{h_\a +2}+t_\b\mult \z_\a +\z_{-\b}\mult\z_{\a +\b} \fracds{2h_\a +h_\b +6}{(h_\b +1)(h_\a +h_\b +3)},\\
\nonumber \z_{\a +\b}\mult\z_{-\b}=-t_\a\mult \z_\a \fracds{h_\b -1}{h_\b (h_\a +h_\b +2)}
- t_\b\mult \z_\a \fracds{h_\a +2h_\b +2}{h_\b (h_\a +h_\b +2)}\\
\phantom{\z_{\a +\b}\mult\z_{-\b}=}{}
+ \z_{-\b}\mult\z_{\a +\b} \fracds{(h_\b +2)(h_\b -1)}{h_\b (h_\b +1)} .
\end{gather}

Weight $\beta$:
\begin{gather}
 \z_\b\mult t_\a =t_\a\mult \z_\b -t_\b\mult \z_\b \fracds{1}{h_\b +2}-\z_{-\a}\mult\z_{\a +\b} \fracds{h_\a +2h_\b +6}{(h_\a +1)(h_\a +h_\b +3)} ,\\
 \z_\b\mult t_\b =t_\b\mult \z_\b \fracds{h_\b +4}{h_\b +2}+\z_{-\a}\mult\z_{\a +\b} \fracds{2h_\a +h_\b +6}{(h_\a +1)(h_\a +h_\b +3)} ,\\
 \nonumber \z_{\a +\b}\mult\z_{-\a}=t_\a\mult \z_\b \fracds{2h_\a +h_\b +2}{h_\a (h_\a +h_\b +2)}
+t_\b\mult \z_\b \fracds{h_\a -1}{h_\a (h_\a +h_\b +2)}\\
\phantom{\z_{\a +\b}\mult\z_{-\a}=}{} +  \z_{-\a}\mult\z_{\a +\b} \fracds{(h_\a +2)(h_\a -1)}{h_\a (h_\a +1)}.
\end{gather}

Weight $-\beta$:
\begin{gather}
 t_\a\mult \z_{-\b} =  \z_{-\b}\mult t_\a - \z_{-\b}\mult t_\b   \fracds{1}{h_\b }-\z_{-\a -\b}\mult\z_\a \fracds{h_\a +2h_\b +3}{(h_\a +2)(h_\a +h_\b +2)} ,\\
  t_\b\mult \z_{-\b}   = \z_{-\b}\mult t_\b  \fracds{h_\b+2}{h_\b }+\z_{-\a -\b}\mult\z_\a \fracds{2h_\a +h_\b +6}{(h_\a +2)(h_\a +h_\b +2)} ,\\
 \nonumber \z_\a\mult \z_{-\a -\b}=\z_{-\b}\mult  t_\a \fracds{2h_\a +h_\b +2}{(h_\a+1)(h_\a +h_\b +1)}+\z_{-\b}\mult t_\b  \fracds{h_\a}{(h_\a +1)(h_\a +h_\b +1)}\\
\phantom{\z_\a\mult \z_{-\a -\b}=}{} + \z_{-\a -\b}\mult\z_\a \fracds{h_\a (h_\a  +3)}{(h_\a +1)(h_\a +2)} .
\end{gather}

Weight $-\alpha$:
\begin{gather}
t_\a\mult \z_{-\a} =\z_{-\a}\mult t_\a  \fracds{h_\a+2}{h_\a }-\z_{-\a -\b}\mult\z_\b \fracds{h_\a +2h_\b +6}{(h_\b +2)(h_\a +h_\b +2)} ,\\
 t_\b\mult \z_{-\a} =-\z_{-\a}\mult t_\a \fracds{1}{h_\a }+\z_{-\a}\mult t_\b +\z_{-\a -\b}\mult\z_\b \fracds{2h_\a +h_\b +3}{(h_\b +2)(h_\a +h_\b +2)} ,\\
\nonumber \z_\b\mult\z_{-\a -\b}=-\z_{-\a}\mult t_\a \fracds{h_\b}{(h_\b +1)(h_\a +h_\b +1)}-\z_{-\a}\mult t_\b \fracds{h_\a +2h_\b +2}{(h_\b +1)(h_\a +h_\b +1)}\\
\phantom{\z_\b\mult\z_{-\a -\b}=}{} +  \z_{-\a -\b}\mult\z_\b  \fracds{h_\b(h_\b +3)}{(h_\b +1)(h_\b +2)} .
\end{gather}

Weight $-\alpha-\beta$:
\begin{gather}
 \z_{-\b}\mult\z_{-\a}=-\z_{-\a -\b}\mult t_\a \fracds{1}{h_\a}- \z_{-\a -\b}\mult t_\b \fracds{1}{h_\b}+\z_{-\a}\mult\z_{-\b} ,\\
 \nonumber t_\a\mult\z_{-\a -\b} =\z_{-\a -\b}\mult t_\a \fracds{h_\a h_\b+h_\b^2+h_\a+3h_\b+3}{(h_\b +1)(h_\a +h_\b +1)}
 +  \z_{-\a -\b}\mult t_\b \fracds{h_\b^2+h_\a+4h_\b+3}{h_\b(h_\b +1)(h_\a+h_\b+1)}\\
 \phantom{t_\a\mult\z_{-\a -\b} =}{} -  \z_{-\a}\mult\z_{-\b} \fracds{h_\a +2h_\b +3}{(h_\a +1)(h_\b +1)} ,\\
 \nonumber t_\b\mult\z_{-\a -\b} = \z_{-\a -\b}\mult t_\a \fracds{h_\b-1}{(h_\b +1)(h_\a +h_\b +1)}\\
\nonumber \phantom{t_\b\mult\z_{-\a -\b} =}{} + \z_{-\a -\b}\mult t_\b \fracds{(h_\b -1)(h_\a h_\b+h_\b^2+2h_\a+4h_\b+3)}{h_\b(h_\b +1)(h_\a +h_\b +1)}\\
\phantom{t_\b\mult\z_{-\a -\b} =}{} + \z_{-\a}\mult\z_{-\b} \fracds{2h_\a +h_\b +3}{(h_\a +1)(h_\b +1)} .
\end{gather}

{}Finally, we rewrite the relations of the type 4, that is, of weight~0, in the form of ordering relations. In addition to the
general commutativity relation (subtype (4{\bf a}))
\begin{gather}
 t_\b\mult t_\a =t_\a\mult t_\b ,
\end{gather}
we have three relations of subtype (4{\bf b}):
\begin{gather}
\nonumber \z_\a\mult \z_{-\a} = h_\a -t_\a\mult t_\a \fracds{1}{h_\a}+\z_{-\a}\mult\z_\a \fracds{h_\a (h_\a +3)}{(h_\a +1)(h_\a +2)}
 -  \z_{-\b}\mult\z_\b  \fracds{h_\b +3}{(h_\b +2)(h_\a +h_\b +2)} \\
\phantom{\z_\a\mult \z_{-\a} =}{} + \z_{-\a -\b}\mult\z_{\a +\b} \fracds{h_\a (h_\a +h_\b +4)}{(h_\a +1)(h_\b +1)(h_\a +h_\b +3)} ,\label{zazma-rel}
\\
 \nonumber \z_\b\mult \z_{-\b}=h_\b -t_\b\mult t_\b \fracds{1}{h_\b}-\z_{-\a}\mult\z_\a \fracds{h_\a +3}{(h_\a +2)(h_\a +h_\b +2)}
 + \z_{-\b}\mult\z_\b \fracds{h_\b (h_\b +3)}{(h_\b +1)(h_\b +2)} \\
\phantom{\z_\b\mult \z_{-\b}=}{} +  \z_{-\a -\b}\mult\z_{\a +\b} \fracds{h_\b (h_\a +h_\b +4)}{(h_\a +1)(h_\b +1)(h_\a +h_\b +3)} ,\label{zbzmb-rel}\\
\nonumber \z_{\a +\b}\mult\z_{-\a -\b} =  \fracds{h_\a h_\b (h_\a +h_\b +2)}{(h_\a +1)(h_\b +1)} \\
\phantom{\z_{\a +\b}\mult\z_{-\a -\b} =}{}
 -  \left( t_\a\mult t_\a \fracds{h_\b}{h_\b +1}+2t_\a\mult t_\b +t_\b\mult t_\b \fracds{h_\a}{h_\a +1}
\right) \fracds{1}{h_\a +h_\b +1}\label{frel1} \\
\nonumber \phantom{\z_{\a +\b}\mult\z_{-\a -\b} =}{}  -  \z_{-\a}\mult\z_\a \fracds{h_\a (h_\a +3)}{(h_\a +1)(h_\a +2)(h_\b +1)}
  -  \z_{-\b}\mult\z_\b \fracds{h_\b (h_\b +3)}{(h_\b +1)(h_\b +2)(h_\a +1)}\\
\nonumber \phantom{\z_{\a +\b}\mult\z_{-\a -\b} =}{}
 +  \z_{-\a -\b}\mult\z_{\a +\b} \fracds{h_\a h_\b (h_\a +h_\b +4)
(\th_\a^2\th_\b+\th_\a\th_\b^2+\th_\a^2+\th_\a\th_\b+\th_\b^2)}{(h_\a +1)^2(h_\b +1)^2 (h_\a +h_\b +2)(h_\a +h_\b +3)} ,
\end{gather}
where in one factor in the numerator of the last coef\/f\/icient we returned to the notation $\th_\a=h_\a +1$ and $\th_\b=h_\b +1$ to make
the expression f\/it into the line.

{\bf 2. Relations for $\boldsymbol{{\mathrm{DR}}(\gl_3)}$.} The ordering relations for the reduction algebra ${\mathrm{DR}}(\gl_3)$ are easily restored from
the list \rf{strel1}--\rf{frel1}: the $\gl(3)$ generators $t_1$, $t_2$ and $t_3$, with $t_\a =t_1-t_2$ and $t_\b =t_2-t_3$, can be written as
\begin{gather*}
 t_1=  \fracds{1}{3}(2t_\a +t_\b +I^{(3,t)})  ,\qquad
t_2= \fracds{1}{3}(-t_\a +t_\b +I^{(3,t)}) ,\qquad
t_1= \fracds{1}{3}(-t_\a -2t_\b +I^{(3,t)}) ,
\end{gather*}
where $I^{(3,t)}$ is the image of the central generator of $\gl(3)$, $I^{(3,t)}=t_1+t_2+t_3$. Since $I^{(3,t)}$ is central, one immediately writes relations for
${\mathrm{DR}}(\gl_3)$. We illustrate it on the example of relations between the generator $\z_\a$ and the  $\gl(3)$ generators $t_1$, $t_2$ and $t_3$:
\begin{gather*}
\z_\a t_1=t_1\z_\a\fracds{h_\a +3}{h_\a +2}-t_2\z_\a\fracds{1}{h_\a +2}-\z_{-\b}\z_{\a +\b}
\fracds{h_\b +2}{(h_\b +1)(h_\a +h_\b +3)}  ,\\
 \z_\a t_2=-t_1\z_\a\fracds{1}{h_\a +2}+t_2\z_\a\fracds{h_\a +3}{h_\a +2}+\z_{-\b}\z_{\a +\b}
\fracds{h_\a +h_\b +4}{(h_\b +1)(h_\a +h_\b +3)}  ,\\
 \z_\a t_3=t_3\z_\a-\z_{-\b}\z_{\a +\b}\fracds{h_\a +2}{(h_\b +1)(h_\a +h_\b +3)}  .
\end{gather*}

{\bf 3. Braid group action.} There are two braid group generators, $\q_\a$ and $\q_\b$, for the diagonal reduction algebra ${\mathrm{DR}}(\sl_3)$.
Given the action of $\q_\a$, the action of $\q_\b$ on ${\mathrm{DR}}(\sl_3)$ can be reconstructed by using the automorphism $\omega$, see~\rf{not2a},
arising from the outer automorphism of the root system of $\sl_3$, which exchanges the roots $\a$ and $\b$,
\begin{gather*}
 \q_\b =\omega \q_\a \omega^{-1}  .
\end{gather*}
The action of the automorphism $\omega$ on the Cartan subalgebra $\langle h_\a ,h_\b\rangle$ of the diagonal Lie algebra~$\sl_3$ and on the generators of the reduction algebra ${\mathrm{DR}}(\sl_3)$ reads
\begin{alignat*}{3}
&h_\a\leftrightarrow h_\b ,\qquad && t_\a\leftrightarrow t_\b ,& \nonumber\\
&\z_\a\leftrightarrow \z_\b ,\qquad &&\z_{-\a}\leftrightarrow \z_{-\b} ,& \\
&\z_{\a +\b}\leftrightarrow -\z_{\a +\b} ,\qquad && \z_{-\a -\b}\leftrightarrow -\z_{-\a -\b} .\nonumber
\end{alignat*}

The action of the braid group generator $\q_\a$ on the Cartan subalgebra $\langle h_\a ,h_\b\rangle$ of the diagonal Lie algebra
$\sl_3$ reads:
\begin{gather} \label{ssnaq1}
\q_\a (h_\a )=-h_\a -2 ,\qquad \q_\a (h_\b )=h_\a +h_\b +1 .
\end{gather}
This action reduces to the standard action of the Weyl group for the shifted generators $\th_\a =h_\a +1$ and $\th_\b=h_\b +1$.

The action of $\q_\a$ on the zero weight generators  $\{ t_\a ,t_\b\}$ of the diagonal reduction algebra ${\mathrm{DR}}(\sl_3)$ is given by:
\begin{gather}\label{ssnaq2}
\q_\a (t_\a )=-t_\a  \fracds{h_\a +2}{h_\a},\qquad
\q_\a (t_\b )  = t_\a  \fracds{h_\a +1}{h_\a}+t_\b   .
\end{gather}

Finally, the action of $\q_\a$ on the rest of the generators is
\begin{alignat}{3}
&\q_\a (\z_{\a} )=-\z_{-\a} \fracds{h_\a +1}{h_\a -1} ,\qquad && \q_\a (\z_{-\a} )=-\z_{\a} ,&\nonumber\\
& \q_\a (\z_{\b} )=\z_{\a +\b} ,\qquad && \q_\a (\z_{\a +\b} )=-\z_\b \fracds{h_\a +1}{h_\a} ,& \label{ssnaq3}\\
& \q_\a (\z_{-\a -\b} )=-\z_{-\b},\qquad && \q_\a (\z_{-\b} )=\z_{-\a -\b} \fracds{h_\a +1}{h_\a} .\nonumber
\end{alignat}

The set of ordering relations \rf{strel1}--\rf{frel1} is covariant with respect to the braid group generated by $\q_\a$ and $\q_\b$. ``Covariant'' means
that the elements of the braid group map a relation to a linear over~$\Uh$ combination of relations. For example, the operator~$\q_\a$, up to multiplicative factors
from~$\Uh$, transforms the relation~\rf{zazma-rel} into itself and permutes the relations~\rf{zbzmb-rel} and~\rf{frel1}. Due to the choice~\rf{stord} of the order,
the set of relations \rf{strel1}--\rf{frel1} is invariant with respect to the anti-involution~$\epsilon$. The set of relations \rf{strel1}--\rf{frel1} is covariant under the
involution~$\omega$ as well.

{\bf 4. Central elements of $\boldsymbol{{\mathrm{DR}}(\sl_3)}$.} The degree 1 and degree 2 (in generators $\z_{ij}$) central elements of the
reduction algebra  ${\mathrm{DR}}(\sl_3)$ are:
\begin{gather*}
 {\mathcal{C}}^{\{ {\mathrm{DR}}(\sl_3),1\} }=t_\a (2h_\a +h_\b +6)+t_\b (h_\a +2h_\b +6) ,\\
 {\mathcal{C}}^{\{ {\mathrm{DR}}(\sl_3),2\} }= \fracds{1}{3}(t_\a\mult t_\a +t_\b\mult t_\b+t_\a\mult t_\b +h_\a^2 +h_\b^2 +h_\a h_\b )\\
\phantom{{\mathcal{C}}^{\{ {\mathrm{DR}}(\sl_3),2\} }=}{}
+  \z_{-\a}\mult\z_\a \fracds{h_\a +3}{h_\a +2} +\z_{-\b}\mult\z_\b \fracds{h_\b +3}{h_\b +2}\\
\phantom{{\mathcal{C}}^{\{ {\mathrm{DR}}(\sl_3),2\} }=}{}  +  \z_{-\a -\b}\mult\z_{\a +\b} \fracds{h_\a +h_\b +4}{h_\a +h_\b +3}
\left( 1+ \fracds{1}{h_\a +1}+ \fracds{1}{h_\b +1}\right) +2(h_\a +h_\b ) .
\end{gather*}
Both Casimir operators, $ {\mathcal{C}}^{\{ {\mathrm{DR}}(\sl_3),1\} }$ and ${\mathcal{C}}^{\{ {\mathrm{DR}}(\sl_3),2\} }$ arise from the quadratic Casimir operator~${\mathcal{C}}^{\{ \sl_3,2\} }$ of the Lie algebra $\sl_3$, whose ordered form is
\begin{gather*}
  {\mathcal{C}}^{\{ \sl_3,2\} }=({\mathcal{E}}_{-\a}{\mathcal{E}}_\a +{\mathcal{E}}_{-\b}{\mathcal{E}}_\b +{\mathcal{E}}_{-\a -\b}{\mathcal{E}}_{\a +\b})
+\fracds{1}{3}({\mathcal{H}}_\a^2+{\mathcal{H}}_\b^2+{\mathcal{H}}_\a {\mathcal{H}}_\b)+{\mathcal{H}}_\a +{\mathcal{H}}_\b .
\end{gather*}
The operator ${\mathcal{C}}^{\{ {\mathrm{DR}}(\sl_3),1\} }$ is the image of ${\mathcal{C}}^{\{ \sl_3,2\} }\otimes 1-1\otimes {\mathcal{C}}^{\{ \sl_3,2\} }$ and
the operator ${\mathcal{C}}^{\{ {\mathrm{DR}}(\sl_3),2\} }$ is the image of ${\mathcal{C}}^{\{ \sl_3,2\} }\otimes 1+1\otimes {\mathcal{C}}^{\{ \sl_3,2\} }$.
We calculate ${\mathcal{C}}^{\{ \sl_3,2\} }\otimes 1+1\otimes {\mathcal{C}}^{\{ \sl_3,2\} }$
and replace the multiplication by the product $\mult$.
Here one needs, in addition to \rf{zemod1}--\rf{zemod5}, the expression for $H_{23}\mult H_{23}$ which is obtained by
applying the involution $\omega$ to  \rf{zemod2} and the equality (in the notation of Section~\ref{sectionDerRel}):
\begin{gather*}
H_{12}\mult H_{23}\equiv  H_{12}H_{23}+E_{21}E_{12} \frac{2}{\th_{12} +1}+E_{32}E_{23} \frac{2}{\th_{23} +1}\nonumber\\
\phantom{H_{12}\mult H_{23}\equiv}{} -  E_{31}E_{13}\left(1+\frac{2}{\th_{12} +1}+\frac{2}{\th_{23} +1}\right)\frac{1}{\th_{13}+1} .
\end{gather*}

The central elements ${\mathcal{C}}^{\{ {\mathrm{DR}}(\sl_3),1\} }$ and ${\mathcal{C}}^{\{ {\mathrm{DR}}(\sl_3),2\} }$ are invariant with respect to the braid group:
\begin{gather*}
 \q_\a \bigl({\mathcal{C}}^{\{ {\mathrm{DR}}(\sl_3),i\} }\bigr)={\mathcal{C}}^{\{ {\mathrm{DR}}(\sl_3),i\} }  ,
\qquad \q_\b \bigl({\mathcal{C}}^{\{ {\mathrm{DR}}(\sl_3),i\} }\bigr)={\mathcal{C}}^{\{ {\mathrm{DR}}(\sl_3),i\} }  ,\qquad i=1,2 .
\end{gather*}
The  central elements ${\mathcal{C}}^{\{ {\mathrm{DR}}(\sl_3),1\} }$ and ${\mathcal{C}}^{\{ {\mathrm{DR}}(\sl_3),2\} }$ are invariant with respect to the anti-involution
$\epsilon$ and the involution $\omega$ as well.

{\bf 5. Diagonal reduction algebra ${\mathrm{DR}}(\sl_2)$.}
For the diagonal reduction algebra of $\sl_2$ we use the following notation:
\begin{gather*}
\z_+:=z_\a ,\qquad z_-:=\z_{-\a} ,\qquad t:=t_\a ,\qquad h:=h_\a.
\end{gather*}
The cut provides the following description of the algebra $\Z_2$ with generators $\z_+ $, $\z_{-}$ and $t$:
\begin{gather}  \z_+\mult t  =t\mult \z_+ \fracds{h +4}{h +2}  ,\label{rsl2pe}\\
\z_+\mult \z_{-} =h -t\mult t  \fracds{1}{h}+\z_{-}\mult\z_+  \fracds{h (h +3)}{(h +1)(h +2)} ,\label{rsl2pm}\\
 t\mult \z_{-}  =\z_{-}\mult t   \fracds{h+2}{h}  .\label{rsl2po}
 \end{gather}

The Casimir operators for ${\mathrm{DR}}(\sl_2)$ are
\begin{gather}\label{cas-sl21}
{\mathcal{C}}^{\{ {\mathrm{DR}}(\sl_2),1\} } :=(h+2)t ,\\
\label{cas-sl22}  {\mathcal{C}}^{\{ {\mathrm{DR}}(\sl_2),2\} }:=\z_-\mult\z_+ \fracds{(h+3)}{(h+2)}+t\mult t  \fracds{1}{4}+\fracds{h(h+4)}{4}.
\end{gather}
Both operators, ${\mathcal{C}}^{\{ {\mathrm{DR}}(\sl_2),1\} }$ and ${\mathcal{C}}^{\{ {\mathrm{DR}}(\sl_2),2\} }$ arise from the quadratic Casimir operator
${\mathcal{C}}^{\{\sl_2,2\} }$ of the Lie algebra $\sl_2$,
\[
{\mathcal{C}}^{\{\sl_2,2\} }={\mathcal{E}}_{-}{\mathcal{E}}_+ +\frac{1}{4}{\mathcal{H}}({\mathcal{H}}+2)  ,
\]
${\mathcal{C}}^{\{ {\mathrm{DR}}(\sl_2),1\} }$ is the image of ${\mathcal{C}}^{\{ \sl_2,2\} }\otimes 1-1\otimes {\mathcal{C}}^{\{ \sl_2,2\} }$ and
${\mathcal{C}}^{\{ {\mathrm{DR}}(\sl_2),2\} }$ is the image of ${\mathcal{C}}^{\{ \sl_2,2\} }\otimes 1+1\otimes {\mathcal{C}}^{\{ \sl_2,2\} }$.

The Casimir operators can be obtained by the cutting also, as explained in Subsection~\ref{subsection3.4}, see Proposition~\ref{proposition2}.
One replaces the  $\sl_3$ generators by the $\gl_3$ generators in the Casimir operators
for $\sl_3$ then cuts and rewrites, using the notation \rf{clih} and \rf{clit}, the result according to the~$\gl_2$ formulas
\begin{gather*}
\begin{split}
& t_1^{(2)}= \fracds{1}{2}\big(t +I^{(2,t)}\big)  ,\qquad t_2^{(2)}=  \fracds{1}{2}\big(-t +I^{(2,t)}\big)  ,\\
& t=t_1^{(2)}-t_2^{(2)} ,\qquad I^{(2,t)}=t_1^{(2)}+t_2^{(2)}  .
\end{split}
\end{gather*}
The cut of ${\mathcal{C}}^{\{ {\mathrm{DR}}(\sl_3),1\} }$ is
\begin{gather}
\frac{3}{2}{\mathcal{C}}^{\{ {\mathrm{DR}}(\sl_2),1\} } + \frac{1}{2}I^{(2,t)}\mult\bigl( I^{(2,h)}+6\bigr)
  -t_3^{(3)}\mult\bigl( I^{(2,h)}+6\bigr)-I^{(2,t)}h_3+2t_3^{(3)}h_3   \label{cuca321}
\end{gather}
and the cut of ${\mathcal{C}}^{\{ {\mathrm{DR}}(\sl_3),2\} }$ is
\begin{gather}
{\mathcal{C}}^{\{ {\mathrm{DR}}(\sl_2),2\} } + \frac{1}{12}I^{(2,t)}\mult I^{(2,t)}+ \frac{1}{12}{I^{(2,h)}}\mult\bigl( I^{(2,h)}+12\bigr)\nonumber\\
\qquad{} - \frac{1}{3}I^{(2,t)}\mult t_3^{(3)}-\left(  \frac{1}{3}I^{(2,h)}+2\right) h_3+\dfrac{1}{3}\bigl(
t_3^{(3)}\mult t_3^{(3)} +h_3^2\bigr)  .  \label{cuca322}
\end{gather}
As expected, the coef\/f\/icients of $(t_3^{(3)})^{\mult i}  h_3^j$ for all $i$ and $j$ in the expressions~(\ref{cuca321}) and~(\ref{cuca322}) are central
elements of the algebra~$\Z_2$.

Due to \rf{ssnaq1},  \rf{ssnaq2} and~\rf{ssnaq3},  the action of the braid group generator reads
\begin{gather}\label{qsl2}   \q (h)=-h-2  ,\qquad
\q (t)=-t  \fracds{h+2}{h}  ,\qquad
\q (\z_{+} )=-\z_{-}  \fracds{h+1}{h-1} ,\qquad
\q  (\z_{-} )=-\z_{+}  .\end{gather}
It preserves the commutation relations of ${\mathrm{DR}}(\sl_2)$. The Casimir operators \rf{cas-sl21} and \rf{cas-sl22} are invariant under the transformation~\rf{qsl2} and under the anti-involution $\epsilon$.

It should be noted that $\q$ can be included in a family of more general automorphisms of the reduction algebra ${\mathrm{DR}}(\sl_2)$.

\begin{lem}
The most general automorphism of the reduction algebra ${\mathrm{DR}}(\sl_2)$ transforming the weights of elements
in the same way as the operator~$\q$ and linear over~$\Uh$ in the genera\-tors~$\z_{+}$,~$\z_{-}$ and $t$ is
\begin{gather}  h\mapsto -h-2  ,\qquad
t\mapsto\beta  t   \frac{h+2}{h} ,\qquad
\z_+\mapsto\z_- \frac{1}{(h-1) \gamma (h)} ,\nonumber\\
\z_-\mapsto   \z_+  (h+3) \gamma (h+2) ,\label{qsl2gen}
\end{gather}
where $\beta =\pm 1$ is a constant and $\gamma (h)$ is an arbitrary function.
\end{lem}

\begin{proof} We are looking for an invertible transformation which preserves the relations \rf{rsl2pe}--\rf{rsl2po} and has the form
\begin{gather}
 h\mapsto f_1(h)  ,\qquad t\mapsto t f_2(h)  ,\qquad \z_+\mapsto\z_-f_3(h)  ,\qquad \z_-\mapsto\z_+f_4(h)\label{tesl2tr}
\end{gather}
with $f_1(h),f_2(h),f_3(h),f_4(h)\in \Uh$. Applying the transformation \rf{tesl2tr} to the relations~\rf{rsl2pe} and~\rf{rsl2po},
we f\/ind (after simplif\/ications) the conditions:
\begin{gather}
(h+2)\bigl( f_1(h)+4\bigr) f_2(h-2)-x\bigl( f_1(h)+2\bigr) f_2(h)=0  ,\label{ttra1}\\
 (h+4)\bigl( f_1(h)+2\bigr) f_2(h)-(x+2) f_1(h) f_2(h+2)=0  .\label{ttra2}
 \end{gather}
Replacing $h$ by $h-2$ in the second equation and then excluding $f_2$ from the system \rf{ttra1}, \rf{ttra2}, we obtain the dif\/ference equation
\begin{gather*} f_1(h)-f_1(h-2)+2=0 ,
\end{gather*}
whose general solution in $\Uh$ is
\begin{gather}
 f_1(h)=-h+c  ,\label{ttra4}
 \end{gather}
where $c$ is a constant.

Applying the transformation \rf{tesl2tr} to the relation \rf{rsl2pm} and collecting the free term and the terms with~$t\mult t$ and $\z_-\mult\z_+$,
we f\/ind (after simplif\/ications)
\begin{gather}
 1+ \fracds{h\bigl( f_1(h)+3\bigr) G(h)}{\bigl( f_1(h)+2\bigr) \bigl( f_1(h)+1\bigr)}=0 ,\label{ttra5}\\
 \fracds{f_2(h)^2}{f_1(h)}+ \fracds{f_1(h)\bigl( f_1(h)+3\bigr) G(h)}{h\bigl( f_1(h)+2\bigr) \bigl( f_1(h)+1\bigr)}=0 ,\label{ttra6}\\
  G(h+2)+ \fracds{h(h+3)f_1(h)\bigl( f_1(h)+3\bigr)}{(h+1)(h+2)\bigl( f_1(h)+2\bigr) \bigl( f_1(h)+1\bigr)} G(h)=0 ,\label{ttra7}
 \end{gather}
where $G(h):=f_3(h)f_4(h-2)$. Excluding $G$ from the system \rf{ttra5}, \rf{ttra6}, we obtain
\begin{gather}
 f_1(h)^2=h^2  f_2(h)^2\qquad {\mathrm{or}}\qquad f_2(h)=\beta  \fracds{f_1(h)}{h}\label{ttra8}
 \end{gather}
with $\beta^2=1$.

The substitution of \rf{ttra4} and \rf{ttra8} into \rf{ttra1} leads to
\begin{gather*} c=-2
\end{gather*}
and it then follows from \rf{ttra5} that
\begin{gather*}
G(h)=\fracds{h+1}{h-1}  .
\end{gather*}
The remaining relation \rf{ttra7} is now automatically satisf\/ied. The proof is f\/inished.
\end{proof}

The Casimir operator ${\mathcal{C}}^{\{ {\mathrm{DR}}(\sl_2),2\} } $ is invariant under the general automorphism~\rf{qsl2gen}.
The Casimir operator ${\mathcal{C}}^{\{ {\mathrm{DR}}(\sl_2),1\} } $ is invariant under the automorphism~\rf{qsl2gen} if\/f $\b =-1$.

The map $\q$ def\/ined by (\ref{qsl2}) is a particular choice of (\ref{qsl2gen}), corresponding to $\beta =-1$ and $\gamma (h)=-\frac{1}{h+1}$.

The map (\ref{qsl2}) is not an involution (but it squares to the identity on the weight zero subspace of the algebra). However, the general map~(\ref{qsl2gen}) squares
to the identity on the whole algebra if\/f the function~$\gamma$ is odd,
\[ \gamma (-h)=-\gamma (h) .
\]

\subsection*{Acknowledgments}

We thank Lo\"{\i}c Poulain d'Andecy for computer calculations for the  algebra~$\Z_4$. We thank Elena Ogievetskaya for the help in preparation
of the manuscript.
The present work was partly done during visits of S.K.\ to CPT and CIRM in Marseille. The authors thank both Institutes. S.K.\ was
supported by the RFBR grant 11-01-00962,  joint CNRS-RFBR grant 09-01-93106-NCNIL,
and by Federal Agency for Science and Innovations of Russian Federation under contract
14.740.11.0347.


\LastPageEnding

\end{document}